\documentclass[11pt,a4paper]{amsart}
\usepackage[all]{xy}
\usepackage{amssymb}

\numberwithin{equation}{section}

\newcommand{\ie}{{\em i.e.}\ }
\newcommand{\cf}{{\em cf.}\ }

\newcommand{\ko}{\: , \;}

\newtheorem{theorem}{Theorem}

\newtheorem{lemma}{Lemma}
\newtheorem{proposition}{Proposition}
\newtheorem{corollary}{Corollary}
\newtheorem{conjecture}{Conjecture}
\newtheorem{question}{Question}

\newcommand{\opname}[1]{\operatorname{\mathsf{#1}}}

\renewcommand{\mod}{\opname{mod}\nolimits}

\newcommand{\ind}{\opname{ind}\nolimits}

\newcommand{\der}{\cd}
\newcommand{\Gr}{\opname{Gr}\nolimits}
\newcommand{\dimv}{\underline{\dim}\,}
\renewcommand{\Im}{\opname{Im}\nolimits}
\newcommand{\inprod}[1]{\langle #1 \rangle}

\newcommand{\cok}{\opname{cok}\nolimits}
\newcommand{\im}{\opname{im}\nolimits}
\renewcommand{\ker}{\opname{ker}\nolimits}
\newcommand{\obj}{\opname{obj}\nolimits}

\newcommand{\Z}{\mathbb{Z}}
\newcommand{\N}{\mathbb{N}}
\newcommand{\Q}{\mathbb{Q}}
\newcommand{\C}{\mathbb{C}}
\newcommand{\fq}{{\mathbb F}_q}
\renewcommand{\P}{\mathbb{P}}

\newcommand{\iso}{\stackrel{_\sim}{\rightarrow}}
\newcommand{\id}{\mathbf{1}}

\newcommand{\Hom}{\opname{Hom}}
\newcommand{\go}{\opname{G_0}}

\newcommand{\Ext}{\opname{Ext}}

\newcommand{\Aut}{\opname{Aut}}
\newcommand{\End}{\opname{End}}
\newcommand{\GL}{\opname{GL}}

\newcommand{\supp}{\opname{supp}}

\newcommand{\ca}{{\mathcal A}}

\newcommand{\cc}{{\mathcal C}}
\newcommand{\cd}{{\mathcal D}}

\newcommand{\cF}{{\mathcal F}}

\newcommand{\ch}{{\mathcal H}}

\newcommand{\ct}{{\mathcal T}}

\newcommand{\cx}{{\mathcal X}}

\newcommand{\wmny}{W_{M,N}^Y}\newcommand{\wnmy}{W_{N,M}^Y}
\newcommand{\eps}{\varepsilon}
\newcommand{\DQ}{{\mathcal D}_Q}
\newcommand{\kos}{K_0^{split}(\cc)}

\newcommand{\tilting}{tilting }
\newcommand{\exceptional}{exceptional }
\renewcommand{\hat}[1]{\widehat{#1}}

\begin{document}

\title{From triangulated categories to cluster algebras\\}
\author{Philippe Caldero}
\address{Institut Camille Jordan, Universit\'e Claude Bernard Lyon I,
69622 Villeurbanne Cedex, France}
\email{caldero@igd.univ-lyon1.fr}
\author{Bernhard Keller}
\address{UFR de math\'ematiques, Universit\'e Denis Diderot -- Paris 7,
2 place Jussieu, 75251 Paris Cedex 05, France}
\email{keller@math.jussieu.fr}

\begin{abstract}
The cluster category is a triangulated category introduced for its
combinatorial similarities with cluster algebras. We prove that a
cluster algebra $\ca$ of finite type can be realized as a Hall
algebra, called \exceptional  Hall algebra, of the cluster
category. This realization provides a natural basis for $\ca$.  We prove
new results and formulate conjectures on `good basis' properties, positivity,
denominator theorems and toric degenerations.
\end{abstract}
\date{version du 09/06/2005}

\maketitle
\section{Introduction}
Cluster algebras were introduced by S.~Fomin and A.~Zelevinsky
\cite{fominzelevinsky1}. They are subrings of the field
$\Q(u_1,\ldots,u_m)$ of rational fractions in $m$ indeterminates,
and defined via a set of generators constructed recursively.
These generators are called {\it cluster variables} and are
grouped into subsets of fixed finite cardinality called {\it
clusters}. The recursion process begins with a pair $({\bf x}, B)$,
called a {\it seed}, where ${\bf x}$ is an initial cluster and $B$ is a
rectangular matrix with integer coefficients. \par

The first aim of the theory was to provide an algebraic framework for
the study of total positivity and of Lusztig/Kashiwara's canonical
bases of quantum groups. The first result is the {\it Laurent
phenomenon} which asserts that the cluster variables, and thus the
cluster algebra they generate, are contained in the Laurent polynomial
ring $\Z[u_1^{\pm 1},\ldots, u_m^{\pm 1}]$.\par

Since its foundation, the theory of cluster algebras has witnessed
intense activity, both, in its own foundations and in its connections
with other research areas. One important aim has been to prove, as in
\cite{BFZ}, \cite{scott}, \ldots\ that many algebras encountered in
the theory of reductive Lie groups have (at least conjecturally) a
structure of cluster algebra with an explicit seed. On the other hand,
a number of recent articles have been devoted to establishing links
with subjects beyond Lie theory.  These links mainly rely on the
combinatorics by which cluster variables are grouped into
clusters. Among the subjects concerned we find Poisson geometry
\cite{GSV}, where clusters are interpreted in terms of integrable
systems, Teichm\"uller theory \cite{FG}, where clusters are viewed as
systems of local coordinates, and tilting theory \cite{BMRRT}
\cite{BMR} \cite{GLS}, where clusters are interpreted as sets of
indecomposable factors of a tilting module. \par

A cluster algebra is said to be of {\it finite type} if the number of
cluster variables is finite. In \cite{fominzelevinsky2}, S.~Fomin
and A.~Zelevinsky
classify the cluster algebras of finite type in terms of Dynkin
diagrams. The cluster variables are then in bijection with the {\it
almost positive} roots of the corresponding root system, \ie the roots
which are positive or opposite to simple roots.  Note that this
classification is analogous to P.~Gabriel's classification of
representation-finite quivers; it is also analogous to the
classification of finite dimensional semisimple Lie algebras in the
theory of Kac-Moody algebras.\par

In finite type, the combinatorics of the clusters are governed by
generalized associahedra. The purpose of the decorated categories of
\cite{MRZ} and, later, of the cluster categories of \cite{BMRRT}
\cite{BMR} \cite{CCS} is to offer a better understanding of these
combinatorics.\par

Let $Q$ be a quiver whose underlying graph is a
simply laced Dynkin diagram and let $\mod kQ$ be the category
of finite-dimensional representations of $Q$ over a field $k$. The
cluster category $\cc$ is the orbit category of the bounded derived category
$\cd^b(\mod kQ)$ under the action of a canonical automorphism.  Thus, it
only depends on the underlying Dynkin diagram and not on the
orientation of the arrows in the quiver $Q$.  With hindsight,
the canonical automorphism is chosen so as to extend the
bijection between indecomposable $kQ$-modules and positive
roots to a bijection between indecomposable objects of
$\cc$ and almost positive roots and, hence, cluster variables.
By the results of \cite{BMRRT}, this bijection also induces
a bijection between clusters and `tilting objects' of $\cc$.
These `coincidences' lead us to the following

\begin{question}\label{question}
Can we realize a cluster algebra of finite type as a `Hall algebra' of
the corresponding category $\cc$?

\end{question}

\noindent
A first result in this direction is the cluster variable formula of
\cite{caldchap}. This formula gives an explicit expression for the
cluster variable associated with a positive root $\alpha$
corresponding to an indecomposable module $M_\alpha$:
The exponents of the Laurent monomials of the cluster variable $X_\alpha$
are provided by the homological form on $\mod kQ$ and the coefficients
are Euler characteristics of Grassmannians of submodules of
$M_\alpha$. In the present article, we use this formula to
provide a more complete answer to question \ref{question}
and to obtain structural results on the cluster algebra, some of
which constitute positive answers to conjectures by Fomin and Zelevinsky.
We first describe our structural results:
\par\noindent * {\bf Canonical basis.} We obtain a $\Z$-basis of the
cluster algebra labelled by the set of so-called \exceptional objects
of the category $\cc$. The results below point to an analogy of this
basis with Lusztig/Kashiwara's dual canonical bases of quantum groups.

\par\noindent * {\bf Positivity conjecture.} We prove that the
Laurent expansion of a cluster variable has positive coefficients, when the seed is associated to a quiver orientation.

\par\noindent * {\bf Good basis property and
Toric degenerations.} We prove that these bases are compatible with
`good filtrations' of the cluster algebra. This provides toric
degenerations of the spectrum of cluster algebras of finite type, in
the spirit of \cite{cald}.

\par\noindent * {\bf Denominator
conjecture.} The formula enables us to prove the following: the
denominator of the cluster variable associated to the positive root
$\alpha$ in its irreducible fraction into polynomials in the $u_i$'s
is $\prod_iu_i^{n_i}$ where $\alpha=\sum_i n_i\alpha_i$ is the
decomposition of $\alpha$ in the basis of simple roots.  \par

Now, the main result of this article is the `cluster multiplication theorem',
Theorem \ref{maintheorem}. This result yields a more complete answer to
Question \ref{question}. It provides a `Hall algebra type'
multiplication formula for the cluster algebra. The main part
    of the paper is devoted to the proof of this formula.\par

Recall that by a result of \cite{keller}, the category
$\cc$ is triangulated. The cluster multiplication formula expresses
the product of two cluster variables associated with objects
$L$ and $N$ of the cluster category in terms of Euler characteristics
of varieties of triangles with end terms $L$ and $N$. Repeated
use of the formula leads to an expression of the structure
constants of the cluster algebra in the basis provided by
the \exceptional objects in terms of Euler
characteristics of varieties defined from the triangles of the
category $\cc$. Thus, the cluster algebra becomes isomorphic
to what we call the `\exceptional Hall algebra' of $\cc$.\par

This theorem can be compared with Peng and Xiao's
theorem \cite{PX}, which realizes Kac-Moody Lie algebras as Hall
algebras of a triangulated quotient of the derived category of a
hereditary category. But a closer look reveals some
differences. Indeed, we use the quotient of the set of triangles
$W_{N,M}^Y$ of the form
$$M\rightarrow Y\rightarrow N\rightarrow M[1],$$ by the automorphism
group $\Aut Y$ of the object $Y$, while Peng and Xiao use the quotient
$W_{N,M}^Y/\Aut(M)\times \Aut(N)$.  Hence, our approach is more a
`dual Hall algebra' approach as in Green's quantum group
realization \cite{Green}.
Another difference is that the associativity of the
multiplication is not proved a priori but results from the
isomorphism with the cluster algebra.

The paper is organized as follows: Generalities and
auxiliary results on triangulated categories and, in particular, on
the cluster category are given in section~2 and in the appendix, where
we prove the constructibility of the sets $W_{N,M}^Y/\Aut(Y)$
described above.

Section~3 deals with the cluster multiplication formula. We first
reduce the proof to the case where the objects involved are
indecomposable. Then the indecomposable case is solved.  Here, the
homology functor from $\cc$ to the hereditary category of quiver
representations plays an essential r\^ole. It allows us to bypass the
`triangulated geometry' of $\cc$, which unfortunately is even out of
reach of the methods of \cite{Toen} \cite{ToenVaquie}, because the
graded morphism spaces of the category $\cc$ are not of finite total
dimension.  The main ingredient of the proof is the Calabi-Yau
property of the cluster category, which asserts a bifunctorial duality
between $\Ext^1(M,N)$ and $\Ext^1(N,M)$ for any objects $M$ and $N$.

In section~4, we use Lusztig's positivity results for
canonical bases \cite{lusztig} \cite{lusztigbook} to prove the
positivity theorem. Then we obtain the denominator theorem.

Section~5 deals with good bases for these cluster algebras. We provide
a basis indexed by \exceptional objects of the category $\cc$,
\ie objects without self-extensions. The cluster variable formula
yields that this basis has a `Groebner basis' behaviour and provides
toric degenerations.

The last part is concerned with conjectures for `non
hereditary' seeds of a cluster algebra of finite type. We formulate a
generalization of the cluster variable formula, and other conjectures
which would follow from it, such as results on positivity, simplicial
fans \ldots\ . We close the article with a positivity conjecture for the
multiplication rule of the \exceptional Hall algebra.

\bigskip
\par\noindent {\bf Acknowledgments.} 
The authors thank F.~Chapoton and C.~Geiss for stimulating discussions.
They are grateful to A.~Zelevinsky for pointing out an error in
an earlier version of this article.

\section{The cluster category}
\begin{subsection}{}\label {modkQ}

Let $\Delta$ be a simply laced Dynkin diagram and $Q$ a quiver
with underlying graph $\Delta$. We denote the set of vertices of
$Q$ by $Q_0$ and the set of arrows by $Q_1$. Let $k$ be a field.
We denote by $kQ$ the path algebra of $Q$ and
by $\mod kQ$ the category of finitely generated right
$kQ$-modules. For $i\in Q_0$, we denote by $P_i$ the
associated indecomposable
projective $kQ$-module and by $S_i$ the associated simple module.
The Grothendieck group $\go(\mod kQ)$ is
free abelian on the classes $[S_i]$, $i\in Q_0$, and is thus isomorphic
to $\Z^n$, where $n$ is the
number of vertices of $Q$. For any object $M$ in $\mod kQ$, the
dimension vector of $M$, denoted by $\dimv(M)$, is the class of
$M$ in $\go(\mod kQ)$.\par

Recall that the category $\mod kQ$
is hereditary, \ie we have $\Ext^2(M,N)=0$ for any objects $M$,
$N$ in $\mod kQ$. For all $M$, $N$ in $\mod kQ$, we put
\[
[M,N]^0=\dim\Hom(M,N),\,[M,N]^1=\dim\Ext^1(M,N),\,\inprod{M,N}=[M,N]^0-[M,N]^1.
\]
In the sequel, for any additive category $\cF$, we denote by
$\ind(\cF)$ the subcategory of $\cF$ formed by a system of
representatives of the isomorphism classes of indecomposable objects
in $\cF$. We know that there exists a partial ordering $\preceq_Q$,
also denoted by $\preceq$, on $\ind(\mod kQ)$ such that
\[
[M,N]^0\not=0\Rightarrow M\preceq N, \hskip 1cm M,\,N\in\ind(\mod kQ).
\]
Denote by $r$ be the cardinality of $\ind(\mod kQ)$. We fix a
numbering $Z_k$, $1\leq k\leq r$, of the objects in $\ind(\mod kQ)$
which is compatible with the ordering.\par
We recall from \cite{Ringel90a} \cite{Ringel90b} that there is a canonical
$\Z$-linear category from which $\mod kQ$ is obtained by
base change.

 \end{subsection}
 \begin{subsection}{}

Denote by $\cd=\DQ=\cd^b(\mod kQ)$ the bounded derived category of
the category of finitely generated $kQ$-modules and by $S$ its
shift functor $M \mapsto M[1]$. As shown in \cite{happel}, the category
$\cd$ is a Krull-Schmidt category, and, up to canonical triangle
equivalence,  it only depends on the underlying graph $\Delta$ of $Q$.
We identify the category $\mod kQ$ with the full subcategory
of $\cd$ formed by the complexes whose homology is concentrated in degree $0$.
We simply call `modules' the objects in this subcategory.
The indecomposable objects of $\cd$ are the
$S^jZ_k$, $j\in\Z$, $1\leq k\leq r$.

We still denote by $\dimv(M)\in\go(\cd)$ the
dimension vector of an object $M$ of $\cd$ in the Grothendieck group
of $\cd$.

Let $\tau$ be the AR-translation of $\cd$. It is the autoequivalence
of $\cd$ characterized by the Auslander-Reiten formula:
\begin{equation}\label{AR}
\Ext_{\cd}^1(N,M)\simeq D\Hom_{\cd}(M,\tau N),
\end{equation}
where $M$, $N$ are any objects in $\cd$ and where $D$ is the functor
which takes a vector space to its dual. The AR-translation $\tau$ is
a triangle equivalence and therefore induces an automorphism of
the Grothendieck group of $\cd$. If we identify this group
with the root lattice of the corresponding root system, the
Auslander-Reiten translation corresponds to the Coxeter
transformation, \cf \cite{BrennerButler76} \cite{Gabriel79}.
\par

We now consider the orbit category $\cd/F$,
where $F$ is the autoequivalence $S\circ\tau^{-1}$. The objects of the category
$\cc=\cc_Q=\cd/F$ are the objects of $\cd$ and the morphisms are
defined by
\[
\Hom_{\cd/F}(M,N)=\coprod_{i\in\Z}\Hom_{\cd}(F^iM,N).
\]
The category $\cc_Q$ was defined in \cite{BMRRT}, \cf also \cite{CCS}
for the A$_n$-case. It is the so-called {\it cluster category}.
Like the derived category $\cd$, up to canonical equivalence,
the cluster category $\cc_Q$ only depends on $\Delta$ and not on the
orientation of the quiver $Q$.

 \end{subsection}
 \begin{subsection}{}

We now review the basic properties of the category $\cc$.
The first two points of the following theorem were proved
in \cite{keller} and and the last two points in \cite{BMRRT}.

 \begin{theorem}\label{basic}
\begin{itemize}
  \item[(i)] The category $\cc$ is triangulated and
  \item[(ii)] the natural functor $\pi$ : $\cd\rightarrow\cc$ is
   a triangle functor.
  \item[(iii)] The category $\cc$ is a Krull-Schmidt category and
  \item[(iv)]  we have $\End_{\cc}(M)=k$ for any indecomposable
              object $M$ of $\cc$.
\end{itemize}
 \end{theorem}
 The shift functor of the triangulated category $\cc$ will still be denoted by
 $S$. Often, we will omit the functor $\pi$ from the notations.
With this convention, the objects in
\[
\ind\mod kQ\cup\{SP_i, \,1\leq i\leq n\}
\]
form a set of representatives for the indecomposables of $\cc$, as
shown in Proposition~1.6 of \cite{BMRRT}.
By formula \ref{AR}, we have, for all objects $M$, $N$ of
$\cd$, that
 $$\Ext^1_{\cd}(M,N)\simeq D\Hom_{\cd}(N,\tau M)\simeq D\Hom_{\cd}(\tau^{-1} N,M)\simeq D\Ext^1(FN,M).$$
 Hence, the category $\cc$ is Calabi-Yau of CY-dimension $1$, which
means that the functor $\Ext^1$ is symmetric in the following sense:
 $$\Ext_{\cc}^1(M,N)\simeq D\Ext_{\cc}^1(N,M).$$
In other words, there is an (almost) canonical non degenerate bifunctorial
pairing
$$\phi\,:\,\Ext_{\cc}^1(M,N)\times\Ext_{\cc}^1(N,M)\rightarrow k.$$

\end{subsection}
\begin{subsection}{}\label{decomposition}

In this section, we study the analogy between the (triangulated) category $\cc$
and the (abelian) category $\mod kQ$. We will see that it can be useful to view $\cc$ as
glued together from copies of $\mod kQ'$, where $Q'$ runs through the set
of orientations of $\Delta$.\par

By the previous section, each object $M$ of $\cc$ can be uniquely decomposed in
the following way:
$$M=M_0\oplus SP_M,$$
where $M_0$ is the image under $\pi$ of a `module' in $\cd$, and where $SP_M$
is the image of the shift of a projective module. We will say that an object $M$
of $\cc$ is a {\em module} if $M=M_0$,  and that $M$ is the {\em shift of a projective module}
if $M=SP_M$.\par

The module $M_0$ can be recovered using the functor
$$H^0=\Hom_{\cc}(kQ_{kQ},?)\,:\,\cc\rightarrow \mod kQ.$$
Indeed, we have $$H^0(M)=H^0(M_0)\oplus H^0(SP_M)=\Hom_{\mod kQ}(kQ_{kQ},M_0)\oplus
\Hom_{\cc}(\oplus_i P_i,SP_M)=M_0,$$
as the last factor is zero. The functor $H^0$ is a homological functor,
\ie it maps triangles in $\cc$ to long exact sequences of $kQ$-modules.

We will deduce the following proposition from Proposition~1.7 of
\cite{BMRRT}.

\begin{proposition}\label{exact}
Let $M$, $N$ be indecomposable $kQ$-modules. Then
\begin{itemize}
\item[(i)] $\Ext^1_{\cc}(M,N)=\Ext^1_{kQ}(M,N)\coprod \Ext^1_{kQ}(N,M)$
and at least one of the two direct factors vanishes.
\item[(ii)] any short exact sequence of $kQ$-modules
$\xymatrix{0\ar[r]&M \ar[r]^i & Y \ar[r]^p & N \ar[r]&0}$ provides a
(unique) triangle $\xymatrix{M \ar[r]^i & Y \ar[r]^p & N \ar[r]& SM}$ in $\cc$,
\item[(iii)] if $M\preceq N$ and if there exists a triangle
$\xymatrix{M \ar[r]^i & Y \ar[r]^p & N \ar[r] & SM}$ in $\cc$,
then $Y$ is also a module and  there exists a short exact sequence of
$kQ$-modules $\xymatrix{0\ar[r]&M \ar[r]^i & Y \ar[r]^p & N \ar[r]&0}$.
Moreover, if this sequence is non split, the modules $M$ and $N$ are non
isomorphic and are not isomorphic to indecomposable factors of $Y$.
\end{itemize}
\end{proposition}
\begin{proof} Point (i) is proved in proposition~1.7 of
\cite{BMRRT}. For points (ii) and (iii), we may assume that $\Ext^1_{kQ}(N,M)$
does not vanish. This implies that there is a non zero morphism from $M$ to $\tau N$
and thus we have $M \preceq N$. If $\Ext^1_{kQ}(N,M)$
was also non zero, we would also have $M \preceq N$, hence $M=N$ and
$\Ext_{kQ}(M,N)=0$, a contradiction. Now it follows from point (i) that the
canonical map
\[
\Ext^1_{kQ}(M,N) \to \Ext^1_{\cc}(M,N)
\]
is bijective. Points (ii) and the first assertion of (iii) therefore follow
from the bijection between elements of $\Ext^1$ and classes of short exact sequences
in $\mod kQ$ respectively triangles in $\cc$. For the last
assertion of (iii), we first note that $M$ is not isomorphic to $N$
because no indecomposable module has selfextensions. Now since the sequence is non
split, the map $Y\to N$ factors through the middle term of the
AR-sequence ending in $N$. Therefore each indecomposable factor
of $Y$ strictly precedes $N$. Dually, $M$ must precede each
indecomposable factor of $Y$.
\end{proof}

The following lemma is useful to understand extensions between indecomposable objects
of $\cc$ when  the situation does not fit into the framework of the proposition above.

\begin{lemma}\label{almsplit}
Let $M\rightarrow Y\rightarrow N\rightarrow SM$ be a non split triangle,
where $M$, $N$ are indecomposable objects of $\cc$. Suppose that there exists no
orientation $Q'$ of $\Delta$ such that $M$, $N$ are simultaneously $kQ'$-modules
with $M\preceq_{Q'}N$ via the embedding $\mod kQ'\rightarrow\cd\rightarrow\cc$.
Then, $N=SM$ and $Y=0$.
\end{lemma}
\begin{proof}
Using the AR-translation, it is always possible to choose an embedding of
$\mod kQ$ in $\cd$ such that $M$ is a projective module. By changing orientations,
we can suppose that $M$ is simple projective, say $P_i$. By the hypothesis of the lemma,
$N$ is not a module, and, as $N$ is indecomposable, $N$ is the shift of a projective $SP_j$.
As the triangle above is non split, the morphism $\eps$ : $N\rightarrow SM$ is non zero.
So there exists a non zero morphism from $P_j$ to $P_i$. The assumption on $P_i$ implies $j=i$.
Hence, $N=SM$ and the morphism $\eps$ is an isomorphism. This forces $Y$ to be zero.
\end{proof}
\end{subsection}

\begin{subsection}{}\label{action}

In this section, we give properties of group actions on triangles in
$\cc$.  Actually, most of them are general facts valid in
Krull-Schmidt triangulated categories.\par Let $Y$ be an object of
$\cc$ and $Y=\coprod_j Y_j $ be a decomposition into isotypical
components.  Suppose that, for each $j$, the object $Y_j$ is the sum
of $n_j$ copies of an indecomposable.  By Theorem \ref{basic}(iv), the
endomorphism algebra of each $Y_j$ is isomorphic to a matrix algebra
over $k$.  Therefore, the radical $R$ of $\End_{\cc}(Y)$ is formed by
the endomorphisms $f$ all of whose components $f_{jj} : Y_j\to Y_j$
vanish. We obtain the decomposition
\[
\End_\cc(Y) = \prod_{j} M_{n_j}(k) \oplus R.
\]
Therefore, the group $\Aut(Y)$ of invertible elements of
$\End_{\cc}(Y)$ is isomorphic to $L(Y)\ltimes U$, where
$L(Y)=\prod_j\GL_{n_j}(k)$ and where $U$ is the unipotent group $1+R$.
We denote by $\wnmy$ the set of triples $(i,p,\eta)$ of morphisms such
that $\xymatrix{M \ar[r]^i & Y \ar[r]^p & N \ar[r]^\eta&SM}$ is a
triangle.  The group $\Aut(Y)$ acts on $\wnmy$ by
$g.(i,p,\eta)=(gi,pg^{-1},\eta)$.  There also exists an action of
$\Aut(M)\times\Aut(N)$ on $\wnmy$ given by
$(g,h).(i,p,\eta)=(ig^{-1},hp,Sg\eta h^{-1})$.  In particular, we can
define an action of the group $k^*$ on $\wnmy$ given by
$\lambda.(i,p,\eta)=(\lambda^{-1}i,p,\lambda\eta)$.
\begin{lemma}\label{unipotent}
For any objects $M$, $Y$, $N$ of $\cc$, the group $\Aut(Y)$ acts on $\wnmy$ with unipotent stabilizers.
\end{lemma}
\begin{proof}
Fix a triple $(i,p,\eta)$ in $\wnmy$ and let $g$ in the stabilizer of $(i,p,\eta)$ for the action of $\Aut(Y)$.
We can draw the following commuting diagram
\[
\xymatrix{M \ar[r]^i\ar[d]^{\parallel} &
    Y \ar[r]^p\ar[d]^g & N \ar[r]^\eta\ar[d]^{\parallel} & SM\ar[d]^{\parallel} \\M \ar[r]^i&
    Y \ar[r]^p & N \ar[r]^\eta& SM.}
\]
As the identity $I$ lies in the stabilizer, we get the solid part of the commutative diagram
 \[
\xymatrix{M \ar[r]^i\ar[d]^{0} &
    Y \ar@{.>}[dl]^f\ar[r]^p\ar[d]^{I-g} & N\ar@{.>}[dl]^h \ar[r]^\eta\ar[d]^{0} & SM\ar[d]^{0} \\M \ar[r]_i&
    Y \ar[r]_p & N \ar[r]_\eta& SM.}
\]
The long exact sequence obtained by applying $\Hom_\cc(Y,?)$ to
the bottom triangle shows that there
is a morphism $f$ with $if=I-g$. Similarly, there is morphism $h$
with $hp=I-g$. Therefore, we have $(I-g)^2=(hp)(if)=0$.
This proves our assertion.
\end{proof}

Consider the third projection $p_3$ : $\wnmy\rightarrow\Ext^1(N,M)$.
Define $\Ext^1(N,M)_Y$ to be the image of $\wnmy$ under $p_3$.
For any $\eta$ in $\Ext^1(M,N)_Y$, the group $\Aut(Y)$ clearly
acts on $p_3^{-1}(\eta)$.
\begin{lemma}\label{transitive}
For any $\eta$ in $\Ext^1(N,M)_Y$, the action of $\Aut(Y)$ on $p_3^{-1}(\eta)$
is transitive and has unipotent stabilizers.
\end{lemma}
\begin{proof}
Let $(i,p,\eta)$ and $(i',p',\eta)$ be in $p_3^{-1}(\eta)$.
Then, by \cite[TR3]{verdier}, there exists a morphism $f$ which
makes the following diagram commutative.
\[
\xymatrix{M \ar[r]^i\ar[d]^{\parallel} &
    Y \ar[r]^p\ar@{.>}[d]^f & N \ar[r]^\eta\ar[d]^{\parallel} & SM\ar[d]^{\parallel} \\M \ar[r]^{i'}&
    Y \ar[r]^{p'}  & N \ar[r]^\eta& SM}
\]
By the 5-Lemma, $f$ is invertible. The last assertion is Lemma \ref{unipotent}.
\end{proof}

\end{subsection}
\begin{subsection}{}

Let $\kos$ be the Grothendieck group of the underlying additive
category of $\cc$, \ie the free abelian group generated by the
isomorphism classes of indecomposable objects of $\cc$.
For each object $Y$ of $\cc$, we still denote by $Y$ the
corresponding element of $\kos$.\par

For a variety $X$, we define $\chi_c(X)$ to be the
Euler-Poincar\'e characteristic of the \'etale cohomology
with proper support of $X$, \ie we have
\[
\chi_c(X) = \sum_{i=0}^\infty (-1)^i \dim H^i_c(X, \overline\Q_l).
\]
Let $M$, $N$, $Y$ be objects of $\cc$.
By section~\ref{append}, the subset $\Ext^1(M,N)_Y)$
of the vector space $\Ext^1(M,N)$ is constructible. Clearly it is
conic. Thus its projectivization $\P\Ext^1(M,N)_Y)$ is a variety
and has a well-defined Euler characteristic
$\chi_c(\P\Ext^1(M,N)_Y)$.\par

For any pair $(Z_j,Z_i)$ of indecomposable object of $\cc$, we call
{\it elementary vector} associated to this pair any element
$Z_i+Z_j-Y_{ij}$ in $\kos$ such that $Y_{ij}$ is the middle term of a
non split triangle $\xymatrix{Z_i \ar[r] & Y_{ij} \ar[r] & Z_j \ar[r]&
SZ_i}$ in $\cc$.  In this case, the number
$c_{ij}=\chi_c(\Ext^1(\P(Z_j,Z_i)_{Y_{ij}})$ is called the {\it
multiplicity number} of the elementary vector. The following
proposition is an analogue of Theorem~3.3 in \cite{caldschi}.\par

\begin{proposition}\label{chi}
Let $M$, $N$, $Y$ be any objects of $\cc$ such that $Y\not=M\oplus N$.
\begin{itemize}
\item[(i)] If $M+ N-Y$ is an elementary vector associated to the pair
$(Z_j,Z_i)$, then we have
$$\chi_c(\P\Ext^1(N,M)_Y)=z_iz_jc_{ij},$$
where $z_j$, resp. $z_i$, is the multiplicity of the indecomposable
component $Z_j$, resp. $Z_i$, in $N$, resp. $M$, and where $c_{ij}$ is
multiplicity number associated to the elementary vector.\par
\item[(ii)] If $M+N-Y$ is not an elementary vector, then
$$\chi_c(\P\Ext^1(N,M)_Y)=0.$$
\end{itemize}
\begin{proof} Let $A$, $B$ be two objects of $\cc$.
Since $\cc$ is triangulated, for each morphism $i$ from
$A$ to $B$, there is a triangle $\xymatrix{A
\ar[r]^i & B \ar[r] & C \ar[r]& SA}$ and the
object $C$ is unique up to (non unique) isomorphism.
Hence, we have a partition
$$\Hom_{\cc}(A,B)=\coprod_C\Hom_{\cc}(A,B)_C,$$
where $C$ runs through the isomorphism classes of $\cc$ and
$$\Hom_{\cc}(A,B)_C=\{i\in\Hom_{\cc}(A,B)\,|\, \mbox{There is a
triangle }\xymatrix{A \ar[r]^i & B \ar[r] & C \ar[r]& SA} \}.$$
Now, by \cite[TR2]{verdier} and Lemma~1.3 of \cite{happel}, we have
$\Ext^1(N,M)_Y=\Hom(N,SM)_{SY}$. This implies that
\begin{equation}\label{union}
\Ext^1(A,B)=\coprod_C\Ext^1(A,B)_C.
\end{equation}\par
Now we prove (i) and (ii). Denote by $M=\oplus_im_iZ_i$ and
$N=\oplus_j n_jZ_j$ the decompositions of $M$ and $N$ into
indecomposable objects. Consider the action of
$G:=\Aut(M)\times\Aut(N)$ on $\P\Ext^1(N,M)$. The group $G$ contains
as a subgroup
$$H=\prod_i\Aut(m_iZ_i)\times\prod_j\Aut(n_jZ_j)=
\prod_i\GL_{m_i}\times\prod_j\GL_{n_j}.$$
Under the isomorphism
$$\Ext^1(N,M)=\oplus_{(i,j)}\Ext^1(Z_j,Z_i)\otimes M_{n_j\times
m_j},$$ the action of $H$ on $\Ext^1(N,M)$ corresponds to the product
of the canonical actions of the $\GL_{m_i}\times\GL_{n_j}$ on matrices. Let
$T$ be the torus of $G$ corresponding to the diagonal matrices. We
have $T\subset H\subset G$ and $T$ acts on $\P\Ext^1(N,M)$ with set of
invariants
\begin{equation}
\P\Ext^1(N,M)^T=\coprod_{(i,j)}\coprod_{r,s}\P\Ext^1(Z_j,Z_i)\otimes k E_{rs},
\end{equation}
where, for fixed $(i,j)$, the $E_{rs}$, $1\leq r\leq n_j$,
$1\leq s\leq m_i$, are the
elementary matrices in $M_{n_j\times m_i}$.\par Set ${\mathcal
X}:=\P\Ext^1(N,M)_Y$. As $G$ acts on $ {\mathcal X}$, we have an
action of the torus $T$ and
$${\mathcal X}^T=\P\Ext^1(N,M)^T\cap{\mathcal X}.$$

It follows from Bia\l ynicki-Birula's results \cite{BialynickiBirula73}
that we have $\chi_c(\cx)=\chi_c(\cx^T)$. We will compute $\chi_c(\cx^T)$.
If ${\mathcal X}^T$ is empty, we have
$\chi_c({\mathcal X})=\chi_c({\mathcal X}^T)=0$.
Let us suppose from now on that ${\mathcal X}^T$ is non empty.
Let the element $[\eta]\otimes E_{rs}$ be in
${\mathcal X}^T$ with $[\eta]\in\P\Ext^1(Z_j,Z_i)$. As $G$, hence $H$,
acts on ${\mathcal X}$, the element $[\eta]\otimes E_{r's'}$ is also in
${\mathcal X}^T$ for any $r'$, $s'$, $1\leq r'\leq n_j$, $1\leq s'\leq
m_i$. Now, by \ref{union}, there exists a unique object $Y_{ij}$
such that $[\eta]$ belongs to $\P\Ext^1(Z_j,Z_i)_{Y_{ij}}$. And so, in $\kos$,
we have
\begin{equation}
Y=M+N-Z_i-Z_j+Y_{ij}.
\end{equation}
This implies that $M+N-Y$ is an elementary vector. Moreover, for any
element $[\eta']$ of $\P\Ext^1(Z_j,Z_i)_{Y_{ij}}$, the
element $[\eta']\otimes E_{kl}$ belongs to $\P\Ext^1(N,M)_{Y'}$, where
$$Y'=M+N-Z_i-Z_j+Y_{ij}=Y.$$
Thus,  $Y'$ is isomorphic to $Y$. So we have proved the inclusion
$$\coprod_{k,l}\P\Ext^1(Z_j,Z_i)_{Y_{ij}}\otimes E_{kl}\subset
{\mathcal X}^T.$$ The proposition will follow if we prove the reverse
inclusion. Suppose that $[\eta']$ belongs to
$\P\Ext^1(Z_j,Z_i)_{Y_{ij}'}$ with $[\eta']\otimes E_{kl}$ belonging to
${\mathcal X}^T$. Then, we obtain $Y=M+N-Z_i-Z_j+Y_{ij}'$, which
implies $Y_{ij}'=Y_{ij}$. Now, in order to finish the proof of the reverse
inclusion, it is sufficient to prove that if we have two triangles
$$\xymatrix{M \ar[r] & Y \ar[r] & N \ar[r]^\eta& SM}, \,\,\xymatrix{M
\ar[r] & Y \ar[r] & N \ar[r]^{\eta'}& SM},$$
with $\eta\in\Ext^1(Z_j,Z_i)$ and $\eta'\in\Ext^1(Z_{j'},Z_{i'})$, then
$(i,j)$ equals $(i',j')$. Let us prove this fact.  The first triangle implies
that $M+N-Y=Z_j+Z_i-Y_{ij}$ for some object $Y_{ij}$.  The second one
implies that $M+N-Y=Z_{j'}+Z_{i'}-Y_{i'j'}$ for some object
$Y_{i'j'}$.  Comparing both equalities in $\kos$ we get
$Z_j+Z_i-Y_{ij}=Z_{j'}+Z_{i'}-Y_{i'j'}$. It follows from
Lemma~\ref{almsplit} and part (iii) of Proposition~\ref{exact} that
the modules $Z_j$ and $Z_i$ are non isomorphic and are not isomorphic
to indecomposable factors of $Y_{ij}$, and similarly for $Z_{i'}$,
$Z_{j'}$ and $Y_{i'j'}$.  Therefore, the signs in the equality imply
that $(i,j)=(i',j')$ or $(i,j)=(j',i')$.  In the second case, we have
$Y_{ij}=Y_{i'j'}$, but it follows from Lemma~\ref{almsplit} and part
(iii) of Proposition~\ref{exact}, that we cannot simultaneously
have two triangles
$$\xymatrix{Z_i \ar[r] & Y_{ij} \ar[r] & Z_j \ar[r]^\eta& SZ_i}, \hbox{ and }\,
\xymatrix{Z_j \ar[r] & Y_{ij} \ar[r] & Z_i \ar[r]^{\eta'}& SZ_j}.$$
Hence, we must have $(i,j)=(i',j')$ and we are done.
\end{proof}

\end{proposition}

\end{subsection}

\section{The cluster multiplication theorem}
\begin{subsection}{}

 We now present the main theorem. The subsections \ref{restriction},
 \ref{mainsection1}, \ref{mainsection2} are devoted to the proof of
 the theorem.\par
 For any $kQ$-module $M$, let  $\Gr(M)$ be the Grassmannian
of $kQ$-submodules of $M$, and let $\Gr_e(M)$ be the $e$-Grassmannian of $M$,
\ie the variety of submodules of $M$ with dimension vector $e$.\par
 Following \cite{caldchap} we define
\[
X_? : \obj(\cc_Q) \to \Q(x_1,\ldots,x_n) \ko M \mapsto X_M
\]
to be the unique map with the following properties:
\begin{itemize}
\item[(i)] $X_M$ only depends on the isomorphism class of $M$,
\item[(ii)] we have
\[
X_{M\oplus N} = X_M X_N
\]
for all $M,N$ of $\cc_Q$,
\item[(iii)] if $M$ is isomorphic to $SP_i$ for the $ith$
indecomposable projective $P_i$, we have
\[
X_M=x_i \ko
\]
\item[(iv)] if $M$ is the image in $\cc_Q$ of an indecomposable
$kQ$-module, we have
\begin{equation}\label{ccformula}
X_M = \sum_{e} \chi_c(\Gr_e(M)) \, x^{\tau(e) - \dimv M + e} \ko
\end{equation}
where $\tau$ is the Auslander-Reiten translation on the Grothendieck group
of $\cd_Q$ and, for $v\in \Z^n$, we put
\[
x^v = \prod_{i=1}^n x_i^{\inprod{\dimv S_i, v}}.
\]
\end{itemize}

From now on, we write $\Ext^1(N,M)$ for $\Ext^1_{\cc}(N,M)$ for any
objects $N$, $M$ of $\cc$ .

\begin{theorem}\label{maintheorem}
For any objects $M$, $N$ of $\cc$, we have
\begin{itemize}
\item[(i)] If $\Ext^1(N,M)=0$, then $X_NX_M=X_{N\oplus M}$,\par
\item[(ii)] If $\Ext^1(N,M)\not=0$, then
\[
\chi_c(\P\Ext^1(N,M)) X_N X_M = \sum_Y (\chi_c(\P\Ext^1(N,M)_Y) +
\chi_c(\P\Ext^1(M,N)_Y)) X_Y \ko
\]
where $Y$ runs through the isoclasses of $\cc$.
\end{itemize}
\end{theorem}
Note that the first part is true by the definition of the map
$M\mapsto X_M$.
\end{subsection}

\begin{subsection}{}\label{restriction}

We prove here that Theorem~\ref{maintheorem} is true if it is true for
{\em indecomposable $kQ$-modules} $M$, $N$, and for any orientation $Q$
of $\Delta$.  So we suppose the following: For all indecomposable
$kQ$-modules $Z_i$, $Z_j$ with $\Ext^1(Z_j,Z_i)\not=0$, we have
\[
\chi_c(\P\Ext^1(Z_j,Z_i)) X_{Z_j} X_{Z_i} = \sum_{Y_{ij}}
(\chi_c(\P\Ext^1(Z_j,Z_i)_{Y_{ij}}) + \chi_c(\P\Ext^1(Z_i,Z_j)_{Y_{ij}}))
X_{Y_{ij}}. \leqno(*)
\]
By Proposition~2.6 of \cite{caldchap}, the formula above is also true
for all indecomposable objects $Z_i$, $Z_j$ such that $Z_i=\tau Z_j=S
Z_j$. Hence, by Lemma \ref{almsplit}, our hypothesis implies that the
theorem is true for all indecomposable objects $M$ and $N$ of
$\cc_Q$.\par Suppose now that $M$ and $N$ are arbitrary objects of
$\cc$, with decompositions into indecomposables $M=\oplus m_iZ_i$,
$N=\oplus n_jZ_j$.  Then we have
\[
\chi_c(\P(\Ext^1(N,M)))=\sum m_in_j\chi_c(\P\Ext^1(Z_j,Z_i)).
\]

Moreover, if an object $Y$ is such that $\chi_c(\P\Ext^1(N,M)_Y) $ is non zero,
then, by Proposition~\ref{chi}, there exist  unique objects $Z_i$, $Z_j$, $Y_{ij}$
such that
\[
\chi_c(\P\Ext^1(N,M)_Y)=m_in_j\chi_c(\P\Ext^1(Z_j,Z_i)_{Y_{ij}})
\]
and $Y_{ij}\oplus (n_j-1)Z_j\oplus (m_i-1)M_i=Y$. Hence, by multiplying all equalities
(*) by $m_i n_jX_{(n_j-1)Z_j\oplus (m_i-1)M_i}$ and adding all these equalities,
we obtain the formula of Theorem~\ref{maintheorem}.
\end{subsection}

\begin{subsection}{}\label{mainsection1}
From now on, in the rest of the proof of theorem 2, we can suppose
that $M$ and $N$ are indecomposable $kQ$-modules.\par

Suppose $X_{\fq}$ is a family of sets indexed by all finite fields
$\fq$ such that the cardinality $\mid X_{\fq}\mid$ is a polynomial
$P_X$ in $q$ with integer coefficients.  In particular, we will
consider the case where $X_{\Z}$ is a variety defined over $\Z$ and,
for any field $k$, $X_k$ is the corresponding variety defined by base
change. In this case, we will set $\chi_1(X):=P_X(1)$. Then it is a
consequence of the Grothendieck trace formula, see \cite{rei},
\cite[Lemma 3.5]{caldchap}, that $\chi_1(X)$ is exactly the Euler
characteristic $\chi_c(X_{\C})$ of $X_{\C}$.  We apply this to prove
the following Lemma.

\begin{lemma}\label{polychar}
For any indecomposable modules $M$, $N$, any module $Y$, and any
dimension vector $e$, we have
$\chi_c(\P\Ext^1(M,N)_Y)=\chi_1(\P\Ext^1(M,N)_Y)$ and
$\chi_c(\Gr_e(M))=\chi_1(\Gr_e(M))$.
\end{lemma}

\begin{proof}
By part (i) of Proposition \ref{exact}, we have
\[
\Ext^1(M,N)_Y = \Ext^1_{kQ}(M,N)_Y \mbox{ or }
\Ext^1(M,N)_Y = \Ext^1_{kQ}(N,M)_Y. 
\]
Therefore, we know from \cite{riedtmann} that this variety
is obtained by base change from a variety defined over $\Z$ and
that the cardinality of its set of $\fq$-points is polynomial
in $q$. The corresponding facts for $\Gr_e(M)$ were shown
in \cite{caldchap}.
\end{proof}

This lemma allows us to work in the remainder of the section in the
case where the base field $k$ is the finite field $\fq$ with $q$
elements.

Let $N,M$ be indecomposable $kQ$-modules such that
\[
0\neq \Ext^1_{kQ}(N,M).
\]
This implies that $\Ext^1_{kQ}(M,N)$ vanishes and hence that
\[
\Ext^1_{\cc_Q}(N,M) = \Ext^1_{kQ}(N,M)\oplus D\Ext^1_{kQ}(M,N) =
\Ext^1_{kQ}(N,M).
\]
We can suppose that we are in this case.

For a module $Y$, the group $\Aut(Y)$ acts naturally on $\Gr(Y)$ and this action stabilizes
the subvarieties $\Gr(Y)$. If $Y$ is any object, we define the action of
$L(Y):=\Aut(H^0(Y))\times\Aut(SP_Y)$ on $\Gr(H^0(Y))$ by letting the
second factor act trivially.\par

If a group $G$ acts on the sets $X$ and $Y$, then $G$ acts on $X\times
Y$ by $g.(x,y)=(g.x,g.y)$.  We define as usual the set of orbits $X\times_G Y$
for this action. By section \ref{action}, we obtain:
\begin{lemma}\label{chifiber}
Let $M$, $N$, $Y$ be objects of $\cc$. Suppose that
$\Ext^1(N,M)_Y\not=0$.  Then we have
$$\chi_1(\wnmy \times_{(L(Y)\times k^*)}
\Gr(H^0(Y))=\chi_1(\P\Ext^1(N,M)_Y\times\Gr(H^0(Y)).$$ If $Y$ is a
module, then we have $$\chi_1(\wnmy\times_{(\Aut(Y)\times k^*)}
\Gr(Y))=\chi_1(\P\Ext^1(N,M)_Y\times\Gr(Y)). $$
\end{lemma}

Note that we have
\[
\chi_c(\P\Ext^1(N,M)) = \dim \Ext^1(N,M).
\]
so that we could \lq simplify\rq\ the left hand side of the main
formula in Theorem \ref{maintheorem}.  Now, with the help of Lemma
\ref{chifiber}, we see that the specialization of the claimed formula
at $x_i=1$, $1\leq i\leq n$, is just the equality of the Euler
characteristics $\chi_1$ of the following varieties
\[
L=\P\Ext^1(N,M) \times \Gr(N)\times \Gr(M)
\]
and
\[
R=\left(\coprod_Y \wnmy\times_{\Aut(Y)\times k^*} \Gr(Y)\right) \coprod
\left(\coprod_Y \wmny \times_{L(Y)\times k^*} \Gr(H^0(Y)) \right).
\]
More precisely, the left hand side of the equality claimed
in the theorem is a sum of terms
\[
\chi_c((\P\Ext^1(N,M) \times \Gr_e(N)) \times \Gr_f(M))\, x^{\tau e
-\dimv N + e + \tau f - \dimv M + f}.
\]
Now the variety
\[
L(e,f)=(\P\Ext^1(N,M) \times \Gr_e(N)) \times \Gr_f(M)
\]
is the union of the subvariety $L_1(e,f)$ consisting of all those
triples $([\eps], N', M')$ such that there is a diagram
of $kQ$-modules (since $\Ext^1_{kQ}(N,M)=\Ext_{\cc_Q}(N,M)$)
\[
\xymatrix{ \eps: & 0 \ar[r] & M \ar[r]         & Y  \ar[r]        & N  \ar[r] & 0 \\
                 & 0 \ar[r] & M' \ar[u] \ar[r] & Y' \ar[u] \ar[r] & N' \ar[u] \ar[r] & 0
           }
\]
and its complement $L_2(e,f)$. By Proposition \ref{prop1} below, we
see that the cardinality of $L_2(e,f)$ and hence of $L_1(e,f)$ over $\fq$
is polynomial in $q$, so the numbers $\chi_1(L_i(e,f))$, $i=1,2$ make 
sense.

Now, the term
\[
\chi_c(\P\Ext^1(N,M) \times \Gr_e(N)) \times \Gr_f(M))\, x^{\tau e
-\dimv N + e + \tau f - \dimv M + f}
\]
is the sum of
\[
\chi_1(L_1(e,f))\, x^{\tau e
-\dimv N + e + \tau f - \dimv M + f}
\]
and
\[
\chi_1(L_2(e,f))\, x^{\tau e
-\dimv N + e + \tau f - \dimv M + f}.
\]
Now we examine the right hand side of the equality of the theorem:
First note that since we have $\Ext^1_{kQ}(N,M) =
\Ext^1_{\cc_Q}(N,M)$, the set
$\Ext^1(N,M)_Y$ is empty if $Y$ does not occur as the middle
term of an extension
\[
0 \to M \to Y \to N \to 0
\]
in the category of {\em modules}. Therefore, we have
\[
\chi_c(\P\Ext^1(N,M)_Y) X_Y = \sum_g \chi_c(\P\Ext^1(N,M)_Y \times \Gr_g(Y))
\, x^{\tau g - \dimv Y + g}.
\]
For non empty $\Ext^1(N,M)_Y$, by Proposition \ref{exact}, one defines a map
\[
\wnmy \times \Gr_g(Y) \to \coprod_{e+f=g} L_1(e,f), \,(i,p,\eps, Y')\mapsto ([\eps],Y'\cap i(M), p(Y'))
\]
This map descends to a map $\wnmy \times_{\Aut(Y)\times k^*} \Gr_g(Y) \to \coprod_{e+f=g} L_1(e,f)$
which is surjective with fibers which are affine spaces. One
obtains the equality of the Euler characteristics and therefore
\[
\chi_c(\P\Ext^1(N,M)_Y) X_Y = \sum_{e,f}\chi_1(L_1(e,f))\, x^{\tau e
-\dimv N + e + \tau f - \dimv M + f}.
\]
It remains to be proved that
\[
\sum_{e,f} \chi_1(L_2(e,f))\,x^{\tau e
-\dimv N + e + \tau f - \dimv M + f} = \sum_Y
\chi_c(\P\Ext^1(M,N)_Y) X_Y.
\]
For this, we need a characterization of the points in $L_2(e,f)$.
Let $([\eps], N', M')$ be a point of
\[
L(e,f) = \P\Ext^1(N,M) \times \Gr_e(N) \times \Gr_f(M).
\]
Let
\[
\phi : \Ext^1_{\cc_Q}(M,N) \times \Ext^1_{\cc_Q}(N,M) \to k
\]
be the (almost) canonical duality pairing. Let
\[
\beta : \Ext^1(M,N') \to \Ext^1(M,N)\oplus \Ext^1(M',N')
\]
be the map whose components are induced by the inclusions
$M'\subset M$ and $N'\subset N$.
\begin{proposition}\label{prop1} The following are equivalent
\begin{itemize}
\item[(i)] $([\eps], N',M')$ belongs to $L_2(e,f)$.
\item[(ii)] $\eps$ is not orthogonal to $\Ext^1(M,N)\cap \Im \beta$.
\item[(iii)] There is an $\eta\in\Ext^1(M,N)$ such that
$\phi(\eta \eps)\neq 0$ and such that, if
\[
\xymatrix{N \ar[r]^i & Y \ar[r]^p & M \ar[r]^\eta & SN}
\]
is a triangle of $\cc_Q$, then there is a diagram of
$kQ$-modules
\[
\xymatrix{N \ar[r]^{H^0(i)} & H^0(Y) \ar[r]^{H^0(p)} & M \\
          N' \ar[u]\ar[r]   & Y' \ar[u]\ar[r]        & {\;\;M'} \ar[u] \ko
          }
\]
where $Y'$ is a submodule of $H^0(Y)$, $N'$ is the preimage of $Y'$ and
$M'$ the image of $Y'$.
\end{itemize}
\end{proposition}

The proposition will be proved in section~\ref{mainsection2}.
We continue the proof of the theorem.
Recall the equality we have to prove:
\[
\sum_{e,f} \chi_1(L_2(e,f))\,x^{\tau e
-\dimv N + e + \tau f - \dimv M + f} = \sum_Y
\chi_c(\P\Ext^1(M,N)_Y) X_Y.
\]
Each object $Y$ of $\cc_Q$ is isomorphic
to the sum of $H^0(Y)$ with a module $SP_Y$ for some
projective $P_Y$. With this notation, the right hand side equals
\[
\sum_{g,Y} \chi_c(\P\Ext^1(M,N)_Y \times \Gr_g(H^0(Y)))\, x^{\tau g - \dimv
H^0(Y) + g}\, X_{SP_Y}.
\]
To prove the equality, we will define a correspondence $C$ between
\[
L_2= \coprod_{e,f} L_2(e,f) \subset \P\Ext^1(N,M)\times \Gr(N)\times
\Gr(M)
\]
and
\[
R_2=\coprod_{Y,g} R_2(Y,g) \ko \mbox{ where } R_2(Y,g)=\wmny \times_{L(Y)\times k^*} \Gr_g(H^0(Y)) .
\]
Namely, the correspondence $C\subset L_2 \times R_2$ is formed by all pairs
consisting of a point $([\eps],M',N')$ of $L_2$ and a point
$(i,p,\eta, Y')$ of $R_2$ such that
$\phi(\eta \eps)\neq 0$, $N'=H^0(i)^{-1}(Y')$, $M'=H^0(p)(Y')$.
Note that in this situation,
we have a diagram
\[
\xymatrix{N \ar[r]^-{H^0(i)} & H^0(Y) \ar[r]^-{H^0(p)} & M \\
          N' \ar[u]\ar[r]   & Y' \ar[u]\ar[r]        & M'. \ar[u] .
          }
\]
We say that a variety $X$ is an {\em extension of affine spaces}
if there is a vector space $V$ and a surjective morphism
$X \to V$ whose fibers are affine spaces of constant dimension.

\begin{proposition} \label{prop2}
\item[a)] The projection $p_1 : C \to L_2$ is surjective and
its fibers are extensions of affine spaces.
\item[b)] The projection $p_2 : C \to R_2$ is surjective and
its fibers are affine spaces of constant dimension.
\item[c)] If the pair formed by $([\eps], M', N')$ and $(i,p,\eta,
Y')$ belongs to $C$, then
\[
x^{\tau e
-\dimv N + e + \tau f - \dimv M + f} = x^{\tau g - \dimv
H^0(Y) + g}\, X_{SP_Y}\ko
\]
where $Y=H^0(Y)\oplus SP_Y$, $P_Y$ is projective,
$e=\dimv M'$, $f=\dimv N'$, $g=\dimv Y'$.
\end{proposition}

This proposition, which will be proved in
section~\ref{proof-of-prop2}, allows us to conclude: Indeed,
the variety $C$ is the disjoint union of the
\[
C_{e,f,Y,g} = p_1^{-1}(L_2(e,f)) \cap p_2^{-1}(R_2(Y,g))
\]
by parts a) and b) of the proposition, for each non empty $C_{e,f,Y,g}$, we have equality of cardinalities on $\fq$
\[
\#(L_2(e,f))= \sum_{Y,g} q^{m_{e,f,Y,g}}\#(C_{e,f,Y,g})
\mbox{ and } \sum_{e,f} q^{n_{e,f,Y,g}}\#(C_{e,f,Y,g}) =\chi_c (R_2(Y,g)),
\]
where $m_{e,f,Y,g}$ and $n_{e,f,Y,g}$ are (non positive) integers.
 Moreover, if $C_{e,f,Y,g}$
is non empty, then we have the equality of part c) of the
proposition. Thus we have
\begin{align*}
\sum_{e,f} \#(L_2(e,f))\,x^{\tau e
-\dimv N + e + \tau f - \dimv M + f}
&= \sum_{e,f, Y,g} q^{m_{e,f,Y,g}}\#(C_{e,f,Y,g}) \,x^{\tau e
-\dimv N + e + \tau f - \dimv M + f} \\
&= \sum_{e,f, Y,g} q^{m_{e,f,Y,g}}\#(C_{e,f,Y,g}) \, x^{\tau g - \dimv
H^0(Y) + g}\, X_{SP_Y} \\
&= \sum_{e,f, Y,g} q^{n_{e,f,Y,g}}\#(C_{e,f,Y,g}) \, x^{\tau g - \dimv
H^0(Y) + g}\, X_{SP_Y}+Z \\
&= \sum_{Y,g} \#(R_2(Y,g)) \, x^{\tau g - \dimv
H^0(Y) + g}\, X_{SP_Y}+Z,
\end{align*}
where $Z$ is a term which vanishes at $q=1$. Hence,
\begin{align*}
\sum_{e,f} \chi_1(L_2(e,f))\,x^{\tau e
-\dimv N + e + \tau f - \dimv M + f}
&=\sum_{Y,g} \chi_1(R_2(Y,g))\, x^{\tau g - \dimv
H^0(Y) + g}\, X_{SP_Y} \\
&=\sum_Y
\chi_c(\P\Ext^1(M,N)_Y) X_Y.
\end{align*}

\end{subsection}
\begin{subsection}{}\label{mainsection2}
This section is devoted to the proof of  Proposition~\ref{prop1}.

Before proving the equivalence of (i) and (ii), we need
some preparation: consider the diagram
\[
\xymatrix{
S^{-1}N \ar@{.>}[r]^\eps & M \ar@{.>}[r]^\eta & SN \\
S^{-1}N' \ar[u]^{S^{-1}i_{N'}} \ar@{.>}[r]_{\eps'} & M'
\ar[u]^{i_{M'}} \ar@{.>}[r]_{\eta'} & SN' \ar[u]_{S i_{N'}}}
\]
Note that (i) holds if and only if there is no $\eps'$ which makes the
left hand square commutative. We formalize this as follows:
The diagram yields two complexes
\[
\xymatrix{
(S^{-1}N,M') \ar[r]^-{\alpha'} & (S^{-1}N,M)\oplus (S^{-1}N',M')
\ar[r]^-{\beta'} & (S^{-1}N',M) \\
(M',SN) & (M,SN)\oplus (M',SN') \ar[l]_-{\alpha} & (M,SN')
\ar[l]_-{\beta}}
\]
where we write $(,)$ for $\Hom_{\cc_Q}(,)$ and where
\[
\alpha' = \left[ \begin{array}{c} (i_{M'})_* \\ (S^{-1} i_{N'})^*
\end{array} \right] \ko
\beta' = \left[ (S^{-1} i_{N'})^*, -(i_M')_*  \right] \ko
\alpha = \left[ (i_{M'})_*, (S i_{N'})^* \right] \ko
\beta = \left[ \begin{array}{c} (S i_{N'})^* \\ -(i_{M'})_* \end{array}
\right] .
\]
The two complexes are in duality via the pairings
\[
(\eta, \eps) \mapsto \phi(\eta\circ \eps) \mbox{ and }
(\eta',\eps') \mapsto \phi'(\eta'\circ \eps')
\]
given by the (almost) canonical forms
\[
\phi : \Hom_{\cc_Q}(S^{-1}N,SN) \to k \mbox{ and }
\phi' : \Hom_{\cc_Q}(S^{-1}N', SN') \to k.
\]
Now let us prove the equivalence of (i) and (ii). Let $p$ denote
the projection
\[
(S^{-1}N,M)\oplus (S^{-1}N',M') \to (S^{-1}N,M) .
\]
Then (i) says that $\eps$ does not belong to $p(\ker \beta')$.
This holds if and only if $\eps$ is not orthogonal to the orthogonal of
$p(\ker \beta')$ in $(M,SN)\oplus (M', SN')$.  Now the orthogonal
of the image of the map
\[
\xymatrix{\ker\beta' \ar[r] & (S^{-1}N,M)\oplus (S^{-1}N',M')
\ar[r]^-p & (S^{-1},M)}
\]
is the kernel of its transpose
\[
\xymatrix{\cok \beta & (M,SN)\oplus (M',SN') \ar[l] & (M,SN)
\ar[l]}
\]
and this is precisely $(M,SN) \cap \Im\beta$. So (i) holds if and only if
$\eps$ is not orthogonal to $(M,SN)\cap \Im\beta$, which means
that (i) and (ii) are equivalent.

Let us prove that (ii) implies (iii). We choose a morphism
$f: M \to SN'$ such that $\beta(f)$ belongs to $\Ext^1(M,N)$
and is not orthogonal to $\eps$. This means that we have
\[
(Si_{N'}) \circ f = \eta \ko f\circ i_{M'}=0 \ko \phi(\eta
\eps)\neq 0.
\]
Now we form triangles on $\eta: M \to SN$ and $0: M' \to SN'$.
Thanks to the fact that $\eta i_{M'}=0$, we have a
morphism of triangles
\[
\xymatrix{N \ar[r]^i & Y \ar[r]^p & M \ar[r]^\eta & SN \\
        N' \ar[u] \ar[r] & N'\oplus M' \ar[u]\ar[r] & M' \ar[u] \ar[r]^0 & SN'. \ar[u]}
\]
By applying $H^0$ we obtain a morphism of long exact sequences
\[
\xymatrix{H^0(S^{-1}M) \ar[rr]^-{H^0(S^{-1}\eta)}
\ar[drr]^{H^0(S^{-1}f)} & &
    N \ar[r]^{H^0 i} & H^0 Y \ar[r]^{H^0 p} & M \ar[drr]^{H^0(f)} \ar[rr]^{H^0(\eta)} &  & H^0(SN) \\
H^0(S^{-1}M') \ar[u] \ar[rr]_-0 & &
    N' \ar[u] \ar[r] & N'\oplus M' \ar[u]\ar[r] & M' \ar[u] \ar[rr]_-0 &  & H^0(SN').
    \ar[u]}
\]
Now we let $Y'$ be the image of $N'\oplus M'$ in $H^0 Y$. An easy
diagram chase then shows that $N'$ is the preimage of $Y'$ under
$H^0 i$ and $M'$ is the image of $Y'$ under $H^0 p$.

Let us prove that (iii) implies (ii). We are given $\eta$ such
that $\phi(\eta \eps)\neq 0$, a triangle of $\cc_Q$
\[
\xymatrix{N \ar[r]^i & Y \ar[r]^p & M \ar[r]^\eta & SN}
\]
and a diagram of $kQ$-modules
\[
\xymatrix{N \ar[r]^{H^0(i)} & H^0(Y) \ar[r]^{H^0(p)} & M \\
          N' \ar[u]\ar[r]   & Y' \ar[u]\ar[r]        & {\;\;M'} \ar[u] \ko
          }
\]
where $Y'$ is a submodule of $H^0(Y)$, $N'$ is the preimage of $Y'$ and
$M'$ the image of $Y'$. We will show that $\eta$ belongs to
$\Ext^1(M,N)\cap \Im \beta$. For this we consider the larger
diagram
\[
\xymatrix{H^0(S^{-1}M) \ar[rr]^-{H^0(S^{-1}\eta)}
\ar@{.>}[drr] & &
    N \ar[r]^{H^0 i} & H^0 Y \ar[r]^{H^0 p} & M  \ar[rr]^{H^0(\eta)} &  & H^0(SN) \\
H^0(S^{-1}M') \ar[u]  & &
    N' \ar[u] \ar[r] & Y' \ar[u]\ar[r] & M' \ar[u] \ar@{.>}[rru]_-0 &  & H^0(SN').
    \ar[u]}
\]
Here $H^0(S^{-1}\eta)$ factors through $N'$ since its image is
contained in the kernel of $H^0 i$, which is contained in $N'$.
Moreover, $H^0 \eta$ vanishes on $M'$ since $M'$ is contained
in the image of $H^0 p$. Now recall that $M$ and $N$ are indecomposable
$kQ$-modules and $\Ext^1_{kQ}(N,M)\neq 0$. Therefore $M$ is non
injective and so $S^{-1}M = \tau^{-1}M$ is still a module (and not
just an object of $\cc_Q$) and similarly, $N$ is non projective
and so $SN = \tau N$ is still a module. Moreover, $M'$ cannot have
an injective direct factor (since that would also be an injective
direct factor of $M$) and so $S^{-1}M'=\tau^{-1}M'$ is still a
module and similarly $SN'=\tau N'$ is still a module.

We would like to show that $\eta \in \Hom_{\cc_Q}(M,\tau N)$
comes from a morphism of modules. For this, recall that we
have
\[
\Hom_{\cc_Q}(U,V) = \sum_{n\in\Z} \Hom_{\der^b(kQ)}(F^n U,V)
\]
for arbitrary modules $U,V$, where $F=\tau^{-1}S$. Moreover,
if $U$ and $V$ are indecomposable and either $U$ or $V$ does
not lie on a cycle of $\der^b(kQ)$, then by part b) of proposition
1.5 of \cite{BMRRT}, there is at most
one $n\in\Z$ such that
\[
\Hom_{\der^b(kQ)}(F^n U,V)\neq 0.
\]
We know that indecomposable postprojective $kQ$-modules
do not lie on cycles of $\der^b(kQ)$. Thus if $U$ and $V$ are
indecomposable and one of them is postprojective, we have
\[
\Hom_{kQ}(U,V)\neq 0 \Rightarrow \Hom_{kQ}(U,V) \iso \Hom_{\cc_Q}(U,V).
\]
In the case where $Q$ is a Dynkin quiver (which we assume),
all modules are postprojective.

Now $M$ and $\tau N$ are indecomposable and we have
\[
\Hom_{kQ}(M,\tau N) = D\Ext^1_{kQ}(N,M)\neq 0
\]
and hence
\[
\Hom_{kQ}(M,\tau N) = \Hom_{\cc_Q}(M,\tau N).
\]
In particular, $\eta$ comes from a morphism of modules
so that we have $\eta=H^0 \eta$. Since $H^0 \eta$ vanishes
on $M'\subset M$, the composition of $\eta$ with $i_{M'}$
vanishes. It remains to be shown that $\eta$ factors
through $SN'$. Now $\tau^{-1}M$  and $N$ are also indecomposable
and we have
\[
\Hom_{kQ}(\tau^{-1} M, N) = D\Ext^1_{kQ}(N,M)\neq 0
\]
and hence
\[
\Hom_{kQ}(\tau^{-1} M, N) = \Hom_{\cc_Q}(\tau^{-1} M,N).
\]
Thus $\tau^{-1}\eta$ comes from a morphism of modules and
$\tau^{-1}\eta = H^0(\tau^{-1} \eta)$. Since $H^0(\tau^{-1} \eta)$
factors as $i_{N'} f$ for a morphism $f: \tau^{-1} M \to N'$,
we have $\tau^{-1}(\eta) = i_{N'} f$ and therefore
\[
\eta = \tau(i_{N'}) \tau f = (S i_{N'}) (Sf).
\]

\end{subsection}

\begin{subsection}{}\label{proof-of-prop2}

We now give a proof of  Proposition~\ref{prop2}.\par
We prove part a). The projection $p_1 : C \to L_2$ is surjective
by  the equivalence between (i) and (iii) in
Proposition~\ref{prop1}. Let $([\eps], M', N')$ be in $L_2$ and
pick an element of $p_2 (p_1^{-1}([\eps],
M',N'))$. Recall that it is an $(L(Y)\times k^*)$-orbit and by construction,  for each point $(i,p,\eta,Y')$ of the orbit, we have a diagram
\[
\xymatrix{N \ar[r]^{H^0(i)} & H^0(Y) \ar[r]^{H^0(p)} & M \\
          N' \ar[u]\ar[r]   & Y' \ar[u]\ar[r]        & M' \ar[u] \ko
          }
\]
where $N'$ is the preimage of $Y'$ under $H^0(i)$ and $M'$ the
image of $Y'$ under $H^0(p)$. We have a morphism of short exact
sequences
\[
\xymatrix{ 0 \ar[r] & (H^0 i)(N) \ar[r] & H^0 Y \ar[r] & \Im H^0 p
\ar[r] & 0 \\
           0 \ar[r] & (H^0 i)(N') \ar[u]\ar[r] & Y' \ar[u]\ar[r] &
           M' \ar[u]\ar[r] & 0 .
}
\]
Thus, if we fix the triangle
\[
\xymatrix{N \ar[r]^i & Y \ar[r]^p & M \ar[r]^\eta & SN} \ko
\]
then the possible submodules $Y'$ form an affine space endowed
with a simply transitive action by
\[
\Hom_{kQ}(H^0(p)Y', (H^0 i)(N)/ (H^0 i)(N')) = \Hom_{kQ}(M', N/N') \ko
\]
where we have the last equality because  $H^0(p)(Y')=M'$ and $\ker(H^0 i) \subset
N'\subset N$. Thus the space does not depend on the choice of a point $(i,p,\eta,Y')$ in the orbit.
Moreover, the set of the possible straight lines
\[
[\eta] \in \P(\Ext^1_{\cc_Q}(M,N) \cap \Im \beta) \subset
\P\Ext^1(M,N)
\]
is the complement of the hyperplane defined by $\phi(? \cdot
\eta)=0$ inside
\[
\P(\Ext^1_{\cc_Q}(M,N) \cap \Im \beta).
\]
Thus these $\eta$ also form an affine space. Therefore,
the variety $p_1^{-1}([\eps],M',N')$ is an extension of
affine spaces by Lemma \ref{transitive}.

We prove part b). Let $(i,p,\eta, Y')$ be a point of an orbit in $R_2$ and $([\eta], M', N')$
in $p_1 (p_2^{-1}([\eta],Y'))$. Then, $N':=H^0(i)^{-1}(H^0(Y'))$ and $M':=H^0(p)(H^0(Y'))$ only
depend on the choice of the orbit. Thus the set $p_1 p_2^{-1}([\eta],Y')$ is parametrized
by the $[\eps]$ with $\phi(\eta\eps)\neq 0$. These form an
affine space inside $\P\Ext^1_{\cc_Q}(M,N)$

We prove part c). Consider the diagram of $kQ$-modules
\[
\xymatrix{          & N/N' \ar[r]        & (H^0 Y)/Y' \ar[r]  & M/M' &  \\
H^0(S^{-1}M) \ar[r] & N \ar[u]\ar[r]^{H^0(i)}     & H^0 Y \ar[u]\ar[r]^{H^0(p)} & M
\ar[u] \ar[r] & H^0(SN) \\
                    & N' \ar[u]\ar[r]    & Y' \ar[u]\ar[r] & M' \ar[u] &
}
\]
Since $N'$ is the preimage of $Y'$ under $H^0(i)$, the kernels of
$N' \to Y'$ and $H^0(i)$ are isomorphic. We denote both by $K$.
Dually, since $M'$ is the image of $Y'$ under $H^0(p)$, the
cokernels of $(H^0 Y)/Y' \to M/M'$ and of $H^0(p)$ are isomorphic.
We denote both by $C$. Then in the Grothendieck group of
$kQ$-modules, we have the following equalities
\[
\tau(N') + \tau(M') = \tau(Y') + \tau(K) \mbox{ and }
N/N' + M/M' = (H^0 Y)/Y' + C.
\]
Therefore we have
\[
\tau N' - N + N' + \tau M' - M + M' = \tau(Y') +\tau(K) - (H^0
Y)/Y' - C.
\]
In the notation of the proposition, it remains to be shown that
\[
x^{\dimv \tau(K) -\dimv C} = X_{SP_Y}.
\]
Now in fact, it is easy to see that $X_{SP_i}=x^{\,I_i}$, where $I_i$ is the injective module associated to $i$. Hence, for each projective $kQ$-module $P$,
\[
X_{SP} = x^{\,\nu P},
\]
where $\nu$ is the Nakayama functor.
 So, what we have
to prove is the equality
\[
\tau(K) - C = \nu P_Y
\]
in the Grothendieck group of $kQ$-modules. For this, we
first note that by the triangle
\[
N \to Y \to M \to SN
\]
of $\cc_Q$, the module $K$ is a quotient of
$H^0(S^{-1}M)= H^0(\tau^{-1}M)$. Since $M$ is indecomposable
non injective, $\tau^{-1}M$ is still a module so that
$K$ is a quotient of $\tau^{-1}M$ and $\tau K$ a quotient
of $M$. In particular, $\tau K$ is still a module and
$\tau K = H^0(\tau K)$. Thus it suffices to prove that
\[
H^0(\tau K) - C = \nu P_Y
\]
in the Grothendieck group of $kQ$-modules. For this, we
first note that we have a triangle
\[
\xymatrix{N \ar[r]^{H^0(i)} & H^0 Y \ar[r] & \cok(H^0 i)\oplus SK
\ar[r] & SN }
\]
in $\der^b(kQ)$ and thus in $\cc_Q$. Secondly, we have a split
triangle
\[
\xymatrix{SP_Y \ar[r] & Y \ar[r] & H^0 Y \ar[r]^0 & S^2 P_Y}
\]
in $\cc_Q$; and thirdly, we have the triangle
\[
N \to Y \to M \to SN
\]
in $\cc_Q$. Note that $H^0 i$ is the composition of the morphism
$N \to Y$ with the projection $Y \to H^0 Y$. If we form the
octahedron associated with this composition, the three triangles
we have just mentioned appear among its faces, as well as a
new triangle, namely
\[
\xymatrix{ SP_Y \ar[r] & M \ar[r] & \cok(H^0 i)\oplus SK \ar[r] &
S^2 P_Y}.
\]
Note that $S^2 P_Y= \nu P_Y$ by the Calabi-Yau property. If we
apply $H^*$ to this triangle, we obtain the exact
sequence of $kQ$-modules
\[
\xymatrix{0 \ar[r] & M \ar[r] & \cok(H^0 i)\oplus H^0(\tau K)
\ar[r] & \nu P_Y \ar[r] & H^0(\tau M)}.
\]
Since $M$ is an indecomposable module, $\tau M$ is either an
indecomposable non injective module or zero. The image of
$\nu P_Y \to \tau M= H^0 \tau M$ is injective (as a quotient
of an injective module). Hence it is zero and we get an exact
sequence
\[
\xymatrix{0 \ar[r] & M \ar[r] & \cok(H^0 i)\oplus H^0(\tau K)
\ar[r] & \nu P_Y \ar[r] & 0}.
\]
In the Grothendieck group, this yields
\[
0 = M -\cok(H^0 i) - H^0(\tau K) + \nu P_Y = C - H^0(\tau K) + \nu
P_Y
\]
and this is what we had to prove.
\end{subsection}
\begin{subsection}{}
We give here some examples which illustrate Theorem \ref{maintheorem}.
\par\noindent
{\bf Example 1.}\par\noindent
Suppose that $M$ and $N$ are indecomposable objects such that $\dim\Ext_{\cc}^1(N,M)=1$. As in \cite{BMRRT}, let $B$ and $B^*$ be the unique objects such that there exist non split triangles
$M\rightarrow B\rightarrow N\rightarrow SM$ and   $N\rightarrow B^*\rightarrow M\rightarrow SN$. In this case, we have
$$\Ext^1(N,M)_B=k^*,\;\;\Ext^1(M,N)_{B^*}=k^*.$$
The cluster multiplication theorem then asserts that $X_NX_M=X_B+X_{B^*}$. Note that in this particular case, this formula was a conjecture of \cite{BMRRT} and is since a theorem of \cite{CCS2}, \cite{BMR2}.

\noindent
{\bf Example 2.}\par\noindent
  If $Q$ is the following quiver of type A$_2$:
  \begin{equation}
    \xymatrix{1 & 2\ar[l]}.
 \end{equation}
 Set $M=S_1\oplus S_1$,  $N=S_2\oplus S_2$. If $Y$ is an object such that $\P\Ext^1(N,M)_Y$ is not empty then $Y$ is either $S_1\oplus P_2\oplus S_2$ or $P_2\oplus P_2$ and it is an easy exercise to prove that  the cardinality of $\P\Ext^1(N,M)_Y$ on $\fq$ is respectively $q^2+2q+1$ and $q(q^2-1)$. In a dual way, if $Y$ is an object such that $\P\Ext^1(M,N)_Y$ is not empty then $Y$ is either $S_1\oplus S_2$ or $0$ and the cardinality of $\P\Ext^1(N,M)_Y$ on $\fq$ is respectively $q^2+2q+1$ and $q(q^2-1)$. \par
 The cluster multiplication theorem gives:
 $$X_NX_M=X_{S_1\oplus P_2\oplus S_2}+X_{S_1\oplus S_2}.$$
 Then, applying again the formula -- note that we now in the case of the previous example-- yields:
 $$X_NX_M=X_{2P_2}+2X_{P_2}+1.$$
 Note that this can be easily verified by taking squares  in $X_{S_2}X_{S_1}=X_{P_2}+1$.

\noindent
{\bf Example 3.}\par\noindent
We give an example where the two indecomposable objects $M$ and $N$ are such that their first extension space has dimension 2. \par\noindent
  We consider the following quiver $Q$ of type D$_4$:
  \begin{equation}
    \xymatrix{
      & 3 \ar[d] & \\
      1 \ar[r] & 2 & \ar[l] 4.
    }
  \end{equation}
  Set $\alpha_i=\dimv(S_i)$. Let $R$, $S$, $T$, $U$ be the indecomposable $kQ$-modules with respective dimension vectors $\alpha_2+\alpha_3+\alpha_4$, $\alpha_1+\alpha_2+\alpha_4$, $\alpha_1+\alpha_2+\alpha_3$, $\alpha_1+2\alpha_2+\alpha_3+\alpha_4$. Consider the injective module $M=I_2$ and the simple module $N=S_2$.  If $Y$ is an object such that $\P\Ext^1(M,N)_Y$ is not empty then $Y$ is either $U$, $R\oplus P_1$, $S\oplus P_3$ or $T\oplus P_4$  and the cardinality of $\P\Ext^1(M,N)_Y$ on $\fq$ is respectively $q-2$, 1, 1 and 1. Symmetrically, if $Y$ is an object such that $\P\Ext^1(N,M)_Y$ is not empty then $Y$ is either $SP_2$, $SP_1+S_1$, $SP_3+S_3$ or $SP_4+S_4$  and the cardinality of $\P\Ext^1(N,M)_Y$ on $\fq$ is respectively $q-2$, 1, 1 and 1.
  \par\noindent
  Hence, the cluster multiplication theorem
 gives
 $$2X_NX_M=-X_U+X_{P_1}X_R+X_{P_3}X_S+X_{P_4}X_T-X_{SP_2}+X_{SP_1}X_{S_1}+X_{SP_3}X_{S_3}+X_{SP_4}X_{S_4}.$$ Applying the formula again gives
 $$2X_NX_M=-X_U+3(X_U+1)-X_{SP_2}+3(X_{SP_2}+1),$$ and finally
 $$ X_NX_M=X_U+3+X_{SP_2}.$$

 \end{subsection}
\section{Finite cluster algebras and positivity}\label{positivity}
\begin{subsection}{}

In order to go further, we have to recall some terminology on cluster
algebras.  More precise and complete information can be found in
\cite{fominzelevinsky2}.\par

Let $n\leq m$ be two positive integers. We fix the {\em coefficient field}
$\cF_0=\Q(u_{n+1},\ldots u_m)$ and the {\it ambient field} ${\cF}=
\Q(u_1,\ldots,u_m)$, where the $u_i$'s are indeterminates. Let
${\mathbf x}$ be a free generating set of $\cF$ over $\cF_0$ and let
$\tilde B=(b_{ij})$ be an $m\times n$ integer matrix such that the
submatrix $B=(b_{ij})_{1\leq i,j\leq n}$ is antisymmetric. Such a pair
$({\mathbf x},\tilde B)$ is called {\it a seed}.\par

Let $({\mathbf x},\tilde B)$ be a seed and let $x_j$, $1\leq j\leq n$,
be in ${\mathbf x}$.  We define a new seed as follows. Let $x_j'$ be
the element of $\cF$ defined by the {\it exchange relation}:
\[
x_jx_j'=\prod_{b_{ij}>0}x^{b_{ij}}+\prod_{b_{ij}<0}x^{-b_{ij}},
\]
where, by convention, we have $x_i=u_i$ for $i>n$. Set ${\mathbf
x'}={\mathbf x}\cup\{x_j'\}\backslash \{x_j\}$.  Let $\tilde B'$ be
the $m\times n$ matrix given by
\[
b_{ik}'=\begin{cases}-b_{ik}&\hbox{if } i=j \hbox{ or } k=j\\
b_{ik}+\frac{1}{2}( \,|b_{ij}|\, b_{jk}+b_{ij}\, |b_{jk}|\,) & \hbox{
otherwise.}\cr\end{cases}
\]
Then a result of Fomin and Zelevinsky asserts that $({\mathbf
x'},\tilde B')$ is a seed.  It is called the {\it mutation} of the
seed $({\mathbf x},\tilde B)$ in the direction $x_j$.  We consider all
the seeds obtained by iterated mutations. The free generating sets
occurring in the seeds are called {\it clusters}, and the variables
they contain are called {\it cluster variables}. By definition, the {\it
cluster algebra} $\ca({\mathbf x},\tilde B)$ associated to the seed
$({\mathbf x},\tilde B)$ is the $\Z[u_{n+1},
\ldots, u_m]$-subalgebra of $\cF$ generated by
the set of cluster variables.
The {\it Laurent phenomenon}, see \cite{fominzelevinsky1},
asserts that the cluster variables are Laurent polynomials with integer
coefficients in the $x_i$, $1\leq i\leq m$.  So, we have $\ca({\mathbf
x},B)\subset\Z[x_1^{\pm 1},\ldots,x_m^{\pm 1}]$.

\par Except in
section \ref{toric}, we will be concerned with cluster algebras such
that $n=m$, \ie $\tilde B=B$.  Note that an antisymmetric matrix $B$
defines a quiver $Q=Q_B$ with vertices corresponding to its rows
(or columns) and which has $b_{ij}$
arrows from the vertex $i$ to the vertex $j$ whenever $b_{ij}>0$. The
cluster algebra associated to the seed $({\mathbf x},B)$ will be also
denoted by $\ca(Q)$.\par

An important result of \cite{fominzelevinsky2} asserts that a cluster
algebra is finite, \ie has a finite number of cluster variables, if
and only if there exists a seed associated to a quiver of simply laced
Dynkin type. In this case, the Dynkin type is unique.

Now fix a quiver $Q$ of simply laced Dynkin type. Then, by
\cite{caldchap}, the $\Z$-module generated by the variables $X_M$,
where $M$ runs over the set of objects of $\cc_Q$, is an algebra; it
is the cluster algebra $\ca(Q)$ and the cluster variables are the
$X_M$'s, where $M$ runs through the indecomposable objects of
$\cc_Q$.\par

An object $M$ of $\cc$ is called {\it \exceptional} if it has no
selfextensions, \ie $\Ext^1(M,M)=0$.  An object $T$ of $\cc$ is a
{\it\tilting} object if it is \exceptional, multiplicity free, and has
the following maximality property: if $M$ is an indecomposable object
such that $\Ext^1(M,T)=0$, then $M$ is a direct factor of a direct sum
of copies of $T$.  Note that a \tilting object can be identified with
a maximal set of indecomposable objects $T_1, \ldots, T_n$ such that
$\Ext^1(T_i,T_j)=0$ for all $i,j$.\par

In view of \cite{caldchap}, the main result of \cite{BMRRT} can be
stated as follows: the map $X_?$ : $M\mapsto X_M$ induces a
bijection from the set of \tilting objects to the set of clusters of
$\ca(Q)$.

\end{subsection}
\begin{subsection}{}

Here we prove a positivity theorem that was conjectured in
\cite{fominzelevinsky2}, see also \cite{caldchap}.
\begin{theorem} \label{positivityThm}
For any object $M$ of $\cc_Q$, the variable $X_M$ is in $\Z_{\geq
0}[u_i^{\pm 1}]$.
\end{theorem}
\begin{proof}
It is sufficient to prove that for any $M$ in $\mod kQ$, and for any
$e$ in $\N^n$, we have $\chi_c(\Gr_e(M))\geq 0$. For this, we
recall the construction of the (classical) Hall algebra ${\mathcal H}(Q)$ of
$\mod(Q)$: The algebra ${\mathcal H}(Q)$ is the vector space with basis 
$\{e_M\}_M$, where $M$ runs through the set of isoclasses of 
finite-dimensional $kQ$-modules.
The multiplication rule on $\ch(Q)$ is given by
$$e_M*e_N=\sum_X P_{M,N}(1)e_X,$$
where $P_{M,N}^X$ is the Hall polynomial defined by
\[
P_{M,N}^X(q)=\#\{Y,\,Y \in\Gr (X),\,Y\simeq N,\,X/Y\simeq M\}\mid_{\fq}.
\]
It is known from \cite{Ringel90a} that ${\mathcal H}(Q)$ is an
associative algebra, isomorphic to the enveloping algebra
$U({\mathfrak n})$ of a maximal nilpotent subalgebra ${\mathfrak n}$
of the semisimple Lie algebra ${\mathfrak g}$ associated to the Dynkin
diagram underlying $Q$. Via the isomorphism ${\mathcal H}(Q)\simeq
U({\mathfrak n})$, the basis $\{e_M\}_M$ is identified with a
Poincar\'e-Birkoff-Witt basis of $U({\mathfrak n})$ (in the
sense of \cite{lusztig}).

\par\noindent For any dimension vector $e$, set
$$b_e=\sum_{\dimv(N)=e}e_N\in {\mathcal H}(Q).$$
Then we have
\[
b_{e'}b_e= \sum_{\dimv(N')=e'}e_{N'}\sum_{\dimv(N)=e}e_N=
\sum_{\dimv M=e+e'}\big(\sum_{\dimv N'=e',\,\dimv N=e}P_{N',N}^M(1)\big)e_M.
\]
Hence,
$$b_{e'}b_e=\sum_{\dimv m=e+e'}\chi_c(\Gr_e(M))e_M,$$ by Lemma
\ref{polychar}. Now by \cite[7.3]{lusztig}, for any dimension vector
$e$, the element $b_e$ is in Lusztig's canonical basis of $U({\mathfrak
n})$, when the quantification parameter $q$ is equal to 1. Moreover, by
\cite[par. 14]{lusztigbook}, the product of two elements of the
canonical basis has positive coefficients in its expansion in the
canonical basis. Finally, by \cite[7.11]{lusztig}, an element of the
canonical basis has positive coefficients in its expansion in the
PBW-base $\{e_M\}$. Hence we have $\chi_c(\Gr_e(M))\geq 0$.
\end{proof}
We can have more by noting that the element $b_e$ of the proof is the
element of the canonical basis associated to the dense orbit of the
moduli space of dimension vector $e$. This easily implies that
$\chi_c(\Gr_e(M))>0$ if $\Gr_e(M)\not=\emptyset$. It would be
interesting to prove that the variety $\Gr_e(M)$ has a cellular
decomposition and to find a combinatorial way to calculate its Euler
characteristic.\par\noindent As a particular case of the theorem, we
obtain the
\begin{corollary}
For any quiver $Q$ of simply laced Dynkin type, the cluster variables
of $\ca(Q)$ are Laurent polynomials in the variables $x_i$ with
positive integer coefficients.
\end{corollary}

We can also generalize Fomin and Zelevinky's denominator theorem
\cite{fominzelevinsky2}, see also \cite[Theorem 3.6]{CCS}, to any
quiver $Q$ of simply laced Dynkin type:
\begin{corollary}
Let $M$ be an indecomposable $kQ$-module and set $\dimv M=\sum_i
m_i\,\dimv S_i$. Then the denominator of $X_M$ as an irreducible
fraction of integral polynomials in the variables $u_i$ is $\prod
u_i^{m_i}$.
\end{corollary}
\begin{proof}
By \ref{ccformula} and the positivity theorem, $X_M$ is
a linear combination with
positive integer coefficients of terms $\prod u_i^{n_i}$, $n_i\in\Z$. These
terms are indexed by the set of dimension vectors of submodules
$N$ of $M$, and for each submodule $N$, we have,
by the Serre duality formula, that
$$n_i=-\inprod{N,S_i}-\inprod{S_i,M/N}.$$
So, it is sufficient to prove that
\begin{itemize}
\item[1.] for all $i$, we have $\inprod{N,S_i}+
\inprod{S_i,M/N}\leq(\dimv M)_i$ and
\item[2.] for all $i$, there exists a submodule $N$ such that the
equality holds.
\end{itemize}
First recall that for each module $X$, we have
$\inprod{X,I_i}=\inprod{P_i,X}=(\dimv X)_i$.
Now, as $\mod kQ$ is hereditary, the injective resolution of
$S_i$ yields $$\inprod{N,S_i}\leq \inprod{N,I_i}.$$
Dually, we have $$\inprod{S_i,M/N}\leq \inprod{P_i,M/N}.$$
Adding both inequalities and using the formula above gives the
first point. Now fix a vertex $i$ of the quiver and let $J$
be the set of vertices $j$ such that there exists a path from $i$ to
$j$. Define the subspace $N$ of $M$ to be the sum of the
subspaces $e_jM$, $j\in J$, where $e_j$ is the idempotent
associated with $j$. Then, by construction, $N$ is a submodule of
$M$ with the following
properties: a) $(\dimv N)_j=0$ if
there is a path $j\rightarrow i$, b) $(\dimv M/N)_j=0$ if there is
a path $i\rightarrow j$. Considering the injective resolution
$$0\rightarrow S_i\rightarrow I_i\rightarrow I\rightarrow 0,$$
we obtain the equality $\inprod{N,S_i}=
\inprod{N,I_i}-\inprod{N,I}=\inprod{N,I_i}$, by a).
Dually, property b) implies that $\inprod{S_i,M/M_i}=\inprod{P_i, M/M_i}$.
So we obtain the equality
$\inprod{M_i,S_i}+\inprod{S_i,M/M_i}=(\dimv M)_i$ as required.

\end{proof}
\section{Filtrations and bases}\label{filtration}

As in the previous section, we assume that $Q$ is a quiver
of simply laced Dynkin type. Recall that the elements
of the generating set $X_M$, $M\in\obj(\cc_Q)$, of $\ca(Q)$
can be written
$$X_M=\sum_e\chi_c(\Gr_e(M))\prod x_i^{-\inprod{e, \dimv S_i}-\inprod{\dimv S_i, \dimv M-e}}.$$
We will show that this formula provides `good' filtrations for
finite cluster algebras.
\begin{subsection}{}

Fix a quiver $Q$ of simply laced Dynkin type and let $B$ be the
antisymmetric matrix such that $Q=Q_B$.  We can view $B$ as an
endomorphism of $\go(\mod kQ)$ endowed with the basis $\dimv S_i$,
$1\leq i \leq n$.
\begin{lemma} We have $Be=\sum_i(\inprod{e,\dimv S_i}
-\inprod{\dimv S_i,e})\,\dimv S_i$.
\end{lemma}
\begin{proof} Recall that $B=(b_{ij})$ with
$$b_{ij}=\begin{cases} 1&\hbox{ if } i\rightarrow j \\-1&\hbox{ if }
j\rightarrow i\\ 0 &\hbox { otherwise }\cr\end{cases}.$$ Hence,
$$B(\dimv S_j)=\sum_{i\rightarrow j}\dimv S_i - \sum_{j\rightarrow
i}\dimv S_i.$$ Note now that
$$\inprod{\dimv S_j,\dimv S_i}-\inprod{\dimv S_i,\dimv S_j}=
\begin{cases} 1&\hbox{ if } i\rightarrow j \\-1&\hbox{ if }
j\rightarrow i\\ 0 &\hbox { otherwise. }\cr\end{cases}$$
This proves the lemma.
\end{proof}
We need the
\begin{lemma}
Let $M$ be an indecomposable $kQ$-module and let $N$ be a proper
submodule of $M$. Then, $B(\dimv N)\not=0$.
\end{lemma}
\begin{proof}
Suppose that $N\subset M$ and $\inprod{\dimv N,\dimv
S_i}-\inprod{\dimv S_i,\dimv N}=0$ for all $i$.  This implies that for
any $kQ$-module $X$, we have $\inprod{N,X}-\inprod{X,N}=0$.\par
Suppose first that $N$ has no injective component. Then, $N$ has an
injective hull $I$, with the following property:
$[N,I]^1=0=[I,N]^0$. Hence,
$$\inprod{N,I}-\inprod{N,I}\geq[N,I]^0>0,$$ which contradicts the
above formula.\par

Suppose that $N$ has a non zero injective component $J$. Then since $N$ is
a proper submodule of $M$, the module $J$ is a proper direct
factor of $M$, in contradiction with the assumption that $M$ is
indecomposable.
\end{proof}
For any Laurent polynomial $X$ in the set of variables $\{x_i\}$, the
{\it support} $\supp(X)$ of $X$ is by definition the set of points
$\lambda=(\lambda_1,\ldots,\lambda_n)$ of $\Z^n$ such that the
$\lambda$-component, \ie the coefficient in $\prod_i x^{\lambda_i}$,
of $X$ is non zero.  For any point $\lambda$ in $\Z^n$, identified
with $\go(\mod kQ)$, let $C_\lambda$ be the convex cone with vertex
$\lambda$, and whose edge vectors are generated by $B(\dimv S_i)$.
The previous lemma easily implies the following proposition.
\begin{proposition}\label{cone}
Fix an indecomposable object $M$ of $\cc_Q$, and let $M=M_0\oplus
SP_M$ be its decomposition as in \ref{decomposition}. Then,
$\supp(X_M)$ is in $C_{\lambda_M}$ with $\lambda_M:=(-\inprod{\dimv
S_i,\dimv M_0}+\inprod{\dimv P_M,\dimv S_i})$.  Moreover, the
$\lambda_M$-component of $X_M$ is 1.
\end{proposition}
\end{subsection}
\begin{subsection}{}

The following proposition rephrases a result of \cite{MRZ}.
\begin{proposition}\label{bijection}
The map $\lambda_?$ : $\cc\rightarrow\Z^n$, $M\mapsto \lambda_M$ is
surjective.  Any fiber of a point in $\Z^n$ contains a unique
\exceptional object. The cones generated by the images of \tilting
objects provide a complete simplicial fan.
\end{proposition}
\begin{proof}
We first describe the \exceptional objects of $\cc$. For any $M$ in
$\obj(\cc)$, we denote by $I_M$ the set of $i$ such that $P_i$ is a
component of $P_M$. The following fact is clear :\par The object
$M=M_0\oplus SP_M$ is \exceptional if and only if $M_0$ is
\exceptional and $(\dimv M_0)_i=0$ for any $i$ in $I_M$.\par Recall
now that for any dimension vector $d$ in $\go(\mod kQ)$, there exists
a unique \exceptional module $M_d$ such that $\dimv(M_d)=d$.\par Let
$E$ be the set of \exceptional modules. It decomposes into the
disjoint union $E=\coprod E_I$, where $I$ runs over the set of
partitions of $\{1,\ldots,n\}$ and where $E_I:=\{M\in E,\,I_M=I\}$.
\par For any object $M=M_0\oplus (\oplus_i m_iSP_i)$, we set
$\dimv(M)=\dimv(M_0)-(m_1,\ldots,m_n)$.\par On the one hand, it is
known by \cite{MRZ} that the cones generated by the images under
$\dimv$ of \tilting objects of $\cc$ provide a complete simplicial fan
in $\Z^n$. On the other hand, by the assertion above, $\dimv$ provides
a bijection from $E$ to $\Z^n$, and via this bijection, the map
$\lambda_?$ is piecewise linear -- the domains of linearity are the
$E_I$'s. Moreover, on $E_I$, the matrix of $\lambda$ is triangular and
the diagonal components are
$$d_i=\begin{cases}1&\hbox{ if } i\in I\\-1&\hbox{ if } i\not\in
I\cr\end{cases}.$$ This proves the proposition.
\end{proof}
\end{subsection}
\begin{subsection}{}
 Under the following hypothesis on $Q$,
 we will now define a filtration of the cluster algebra
 $\ca(Q)\subset\Z[x_1^{\pm 1},\ldots,x_n^{\pm 1}]$.
\medskip\par\noindent
HYPOTHESIS: There exists a form $\epsilon$ on $\Z^n$ such that
$$\epsilon(B_Q\dimv S_i)\in \Z_{<0},\;\hbox{for all } i.$$
\par\noindent
Note that for any Dynkin diagram except A$_1$, there exists an
orientation $Q$ satifying our hypothesis.
Indeed, let $Q_{alt}$ be an alternating quiver. Then the
matrix $B_{alt}=(b_{ij})$ associated to this quiver satisfies
$b_{ij}\geq 0$ if $i$ is a source, and $b_{ij}\leq 0$ if $i$ is a
sink. Moreover, each row of $B_{alt}$ is non zero. So we can take
any form $\epsilon$ whose coordinates in the dual basis of $\Z^n$
satisfy $\epsilon_i<0$ if $i$ is a source and $\epsilon_i>0$ if $i$
is a sink. Note also that we have
\[
Be = \sum_i \inprod{e+\tau e, \dimv S_i} \, \dimv S_i
\]
for all $e\in \Z^n$ so that the above hypothesis holds iff the
image of the positive cone of $\go(\mod kQ)$ under the map
$\tau+\id$ is strictly contained in a halfspace.

For any $n$ in $\Z$, set
$$F_n=(\oplus_{\epsilon(\mu)\leq n}\Z \prod x_i^{\mu_i})\cap\ca(Q).$$
Using Proposition \ref{cone} and Proposition
\ref{bijection}, we obtain:
\begin{proposition}
The set $(F_n)_{n\in \Z}$ defines a filtration of $\ca(Q)$.  The
graded algebra asso\-cia\-ted to this filtration is isomorphic to
$\Z[u_1^{\pm 1},\ldots,u_n^{\pm 1}]$.
\end{proposition}
\begin{proof}
As we have $F_nF_m\subset F_{n+m}$, the sequence $(F_n)$ is a filtration of
$\ca(Q)$. Moreover, Proposition \ref{cone} implies that for any
indecomposable module $M$, we have
\begin{equation}\label{gr}
\hbox{gr}X_M=\hbox{gr}\prod_iu_i^{(\lambda_M)_i}.
\end{equation}
The result then follows from Proposition \ref{bijection}.

\end{proof}
This implies
\begin{corollary}\label{base}
For any Dynkin quiver, the set ${\mathcal B}:= \{X_M, \,M\in\obj(\cc),
 \,\Ext^1(M,M)=0\}$ of variables corresponding to \exceptional objets
 of the category $\cc$ is a $\Z$-basis of the cluster algebra
 $\ca(Q)$.
\end{corollary}
\begin{proof}
This is obviously true for a quiver of type A$_1$. Now, for any
quiver, we have $\ca(Q)=\ca(Q')$ for some quiver $Q'$ which
satisfies the hypothesis above. In this case, the proposition together
with formula \ref{gr} imply that  ${\mathcal B}$ is $\Z$-free.\par\noindent
Let us prove now that ${\mathcal B}$ generates $\ca(Q)$ as a $\Z$-module.
We first define a degeneration ordering $\prec_e$ in $\obj(\cc)$. Let
$M$, $M'$ be objects of $\cc$.  We say that there is an {\it
elementary degeneration} $M'\prec_e M$ if $M-M'$ is an elementary
vector in $\kos$. We have

\begin{lemma}\label{elemdeg}

 Let $M$, $M'$ be objects of $\cc$.
\begin{itemize}
\item[a)] We have $M'\prec_e M$ iff there are decompositions
$M=L\oplus U\oplus V$ and $M'=L\oplus E$ where $U$ and $V$ are
indecomposable and $E$ is the middle term of a non split triangle
\[
\xymatrix{U \ar[r] & E \ar[r] & V \ar[r] & SU.}
\]
\item[b)] If we have $M'\prec_e M$, then
 \begin{equation}\label{inequality}
0\leq \dim\Ext^1(M',M') < \dim\Ext^1(M,M).
\end{equation}
\end{itemize}
 \end{lemma}

 \begin{proof}
a) If the condition holds, then $M'-M$ equals the elementary
vector $U+V-E$. Conversely, if $M'-M$ is an elementary vector, we
have $M'-M=U+V-E$ where $U$, $V$ are indecomposable and $E$ is
the middle term of a non split triangle as in the assertion. Then
we have $M'+E = M+U+V$. It follows from Proposition~\ref{exact}
and Lemma~\ref{almsplit} that $U$ and $V$ are not direct factors
of $E$. Moreover, they are non isomorphic since no indecomposable
has selfextensions. Thus, the object $M'$ decomposes as the sum of
some $L$ and $U\oplus V$ so that we obtain $L+U+V+E= M+U+V$ and
$L+E=M$.

b) Let $U$, $V$ and $L$ be as in a).
For any objects $N, N'$ of $\cc$, set $[N,N']^1=\dim\Ext^1_{\cc}(N,N')$.
We view $[?,?]^1$ as a symmetric bilinear form on $\kos$. For any
object $N$ of $\cc$, the long
exact sequence obtained by applying $\Hom(?,SN)$ to the triangle
\[
\xymatrix{U \ar[r] & E \ar[r] & V \ar[r]^e & SU}
\]
shows that we have
\[
[E,N]^1 \leq [U,N]^1+[V,N]^1. \leqno(*)
\]
Moreover, for $N=U$, we have the strict inequality
\[
[E,U] < [U,U]^1+[V,U]^1
\]
since in the sequence
\[
\Hom(V,SU) \to \Hom(E,SU) \to \Hom(U,SU) \ko
\]
the first map is not injective: its kernel contains $e$.
By the inequality $(*)$, we have
\[
[L,L]^1+2[L,E]^1 \leq [L,L]^1+2[L,U+V]
\]
and it only remains to be shown that
\[
[E,E]^1 < [U+V,U+V]^1.
\]
For this, we note that by the above inequalities, we have
\[
[E,E]^1 \leq [E,U]^1+[E,V]^1 \ko
[E,V]^1 \leq [U,V]^1 + [V,V]^1
\mbox{ and }
[E,U]^1< [U,U]^1+[V,U]^1.
\]
\end{proof}

Let us finish the proof of the corollary. It remains to be proved that each
$X_M$ is in $\Z{\mathcal B}$. If $M$ is indecomposable, then $M$ is
\exceptional and hence is in ${\mathcal B}$. If $M$ is not
indecomposable, say $M=M'\oplus M''$, then by the cluster
multiplication theorem and the lemma above, $X_M$ expands into a
$\Q$-linear combination of terms $X_Y$ for objects $Y$ such that $0\leq
\Ext^1(Y,Y)<\Ext^1(M,M)$. By induction, we obtain that $X_M$ is in
$\Q{\mathcal B}$. As the coefficients $\chi_c(\Gr_e(M))$ in the
cluster variable formula are integers, we obtain, using induction
on the filtration, that $X_M$ is in $\Z{\mathcal B}$.

\end{proof}

\end{subsection}
\begin{subsection}{}\label{toric}
In this section, we consider the more general case of finite cluster
algebras associated to a rectangular $m\times n$ matrix $\tilde B$. We
want to prove that our construction provides a toric degeneration of
the spectrum of finite cluster algebras. \par Let $B$ be the
antisymmetric submatrix associated to $\tilde B$. We have a projection
$\pi$: $\ca(\tilde B)\rightarrow\ca(B)$ such that $u_i\mapsto u_i$, if
$1\leq i\leq n$, and $u_i\mapsto 1$, if $n+1\leq i\leq m$. This
projection gives a one-to-one correspondence between the cluster variables
of the two cluster algebras. For any indecomposable object $M$ of $\cc$,
we denote by $\tilde X_M$ the cluster variable such that $\pi(\tilde
X_M)=X_M$. We fix a quiver $Q$ as in the previous section and we
suppose without loss of generality that $B=B_Q$. Let $F_n$ be the
filtration of $\ca(B)$ constructed from $Q$. We now consider the
filtration $\tilde F_n=\pi^{-1}(F_n)$, $n\in\Z$, induced by $\pi$ from
the filtration $F_n$.
\begin{theorem}
The graded algebra gr$\ca(\tilde B)$ associated to the filtration
$\tilde F_n$ is isomorphic to a subalgebra of $\Z[u_i^{\pm 1},\,1\leq
i\leq m]$ generated by a finite set of unitary monomials.
\end{theorem}
\begin{proof}
For any $z$ in $\ca(\tilde B)$, we denote by gr$z$ the corresponding
element in the graded algebra gr$\ca(\tilde B)$. It is sufficient to
prove that for any indecomposable object $M$ of $\cc$, gr$\tilde X_M$
is a unitary monomial in the gr$u_i$'s. This is true for $M=SP_i$ as
in this case, gr$\tilde X_M$ is the monomial $u_i$. Now, we make an
induction with the help of the Hom-ordering $\preceq$ in $\ind\mod kQ$. By
the exchange relation as in \cite[3.4]{caldchap}, we have
$$\tilde X_{\tau(M)}\tilde X_M=p\prod_i\tilde X_{B_i}+q,$$ where $p$,
$q$ are unitary monomials in the $u_i's$, $n+1\leq i\leq m$. In this
relation, $\tau(M)$ and $B_i$ are indecomposable objects which verify
the induction hypothesis. Suppose that gr$\prod \tilde X_{B_i}$ and 1
have the same degree in gr$\ca(\tilde B)$.  Then, the monomial
gr$\prod X_{B_i}$ and $1$ have the same degree in gr$\ca(B)$. But this
would imply that the coefficient of the monomial gr$X_M$ is not $1$,
in contradiction with Proposition \ref{cone}. Hence,
gr$\prod \tilde X_{B_i}$ and 1 do not have the same degree in
gr$\ca(\tilde B)$. This implies that either gr$\tilde
X_{\tau(M)}$gr$\tilde X_M=$gr$ p$ gr$\prod\tilde X_{B_i}$ or gr$\tilde
X_{\tau(M)}$gr$\tilde X_M=$gr$ q$. In both cases, the induction
process is proved.

\end{proof}

This is a classical corollary of the theorem, see \cite{cald}.
\begin{corollary}
The spectrum of a finite cluster algebra has a toric degeneration.
\end{corollary}
\end{subsection}
\begin{subsection}{}
With the help of the basis of Corollary \ref{base}, we can reformulate
Theorem \ref{maintheorem}.  Actually, we can give a complete
realization of the cluster algebra $\ca(Q)$ from the cluster category
$\cc_Q$.\par

We recall the degeneration ordering $\prec_e$ in $\obj(\cc)$ defined
for the proof of Lemma \ref{elemdeg}.  Suppose that there is an
elementary degeneration $M'\preceq_e M$ with elementary vector
$Z_i+Z_j-Y_{ij}$.  Then, by the lemma \ref{elemdeg}, we can define the
ratio
 $$r(M,M'):=cz_iz_j/\dim\Ext^1(M,M),$$ where $c$ is the multiplicity
number associated to the elementary vector, and where $z_i$,
resp. $z_j$, are the multiplicity of $Z_i$, resp. $Z_j$, in $M$.  Let
$\prec$ be the ordering generated by $\prec_e$. \par Note the
surprising fact that in the category $\cc$, the composition of
elementary degenerations can again be an elementary degeneration.\par
Note that, by \ref{inequality}, any chain of elementary
degenerations descending from an object
$M$ is finite.\par For any pair of objects $M$, $N$ of $\cc$ and for
any chain from $N$ to $M$
$$\Gamma\,:\, N=M_0\prec_e M_1\prec_e\ldots\prec_e M_k=M,$$
we set
$$r(M,N,\Gamma)=\prod_{i=1}^k r(M_i,M_{i-1}),\;\; r(M,N)=\sum_\Gamma
r(M,N,\Gamma),$$ where $\Gamma$ runs through the chains from $N$ to
$M$.\par

We define the {\em \exceptional Hall algebra of $\cc$} to be the
rational vector space $\ch_{exc}(\cc)$ of $\Q$-valued functions
on the isomorphism classes of \exceptional objects of $\cc$
endowed with the multiplication given by
$$\chi_M*\chi_N=\sum_{K}r(M\oplus N,K)\,\chi_K \ko$$
where the sum runs over the isomorphism classes of \exceptional
objects $K$ and $\chi_K$ denotes the characteristic function.
From Theorem~\ref{maintheorem} and section~\ref{restriction},
we easily deduce the
\begin{theorem}\label{maintheorem2}
The \exceptional Hall algebra $\ch_{exc}(\cc_Q)$ is an associative
$\Q$-algebra with unit element $\chi_0$.
The map $\chi_M \mapsto X_M$ induces an isomorphism between
the \exceptional Hall algebra $\ch_{exc}(\cc_Q)$ and the cluster
algebra $\ca_Q$.
\end{theorem}

\end{subsection}

\section{Conjectures}
\begin{subsection}{}

The results in \cite{caldchap} and in this article are concerned with
finite cluster algebras, with a fixed seed corresponding to a Dynkin
quiver. We conjecture some generalizations for any seed.

Fix a quiver $Q$ of Dynkin type and a \tilting object $T$ of
$\cc_Q$. Consider the so-called {\it tilted algebra}
$A_T:=\End_{\cc}(T)^{opp}$ and the category $\mod A_T$ of finite
dimensional $A_T$-modules.

We consider the form
$$\inprod{N,M}=\dim\Hom(N,M)-\dim\Ext^1(N,M),\;\;N,M\in \mod A_T.$$
Remark that in general this form does not descend to the Grothendieck
group of the category $\mod A_T$. One defines the antisymmetrized form:
$$\inprod{N,M}_a=\inprod{N,M}-\inprod{M,N},\;\;N,M\in \mod A_T.$$

We know that there exists a seed $({\mathbf x}_T,B_T)$ of the cluster
algebra $\ca(Q)$ associated to the \tilting object $T$, \cf
\cite{BMRRT}. Moreover, by \cite{BMR}, the set $\ind(\mod A_T)$ is in
bijection with the set of cluster variables which do not belong to
${\mathbf x}_T$.
\begin{conjecture}
The form $\inprod{\,,\,}_a$ descends to the Grothendieck group $\go(\mod A_T)$.
Its matrix for the basis $(\dimv S_i)$ is $B_T$.
\end{conjecture}

Set ${\mathbf x}_T=\{x_1,\ldots,x_n\}$.
The following conjecture describes the bijection explained above.
It can be seen as a generalization of Theorem~3.4 of \cite{caldchap}.

\begin{conjecture}\label{ccg}
To any indecomposable module $M$ in $\mod A_T$, we assign
$$X_M:=\sum_e\chi_c(\Gr_e(M))\prod_i x_i^{(B_Te)_i-\inprod{S_i,M}}.$$
Then the set $\{X_M,\,M\in \ind\mod A_T\}$ is exactly the set of cluster
variables which do not belong to ${\mathbf x}_T$.
\end{conjecture}
Via a conjectural extension of theorem~\ref{positivityThm}, 
this conjecture would yield
positivity properties.

\end{subsection}

\begin{subsection}{}

As before, fix a \tilting    object $T=\oplus_{i=1}^n T_i$ of $\cc_Q$, with $T_i$ indecomposable.
By the discussion above, the set $\ind(\cc)$ can be seen as a disjoint union
$$\ind(\cc)=\ind\mod A_T\coprod\{T_i,\,1\leq i\leq n\}.$$
Hence, as in \ref{decomposition}, each object $M$ of $\cc$ has a unique decomposition
$M=M_0\oplus T_M$, where $M_0$ is in $\mod A_T$ and where $T_M$ is
a direct factor of a sum of copies of $T$.

Suppose that Conjecture \ref{ccg} is true. Then, for any
object $M=M_0\oplus (\oplus_im_iT_i)$, we can define the variable

$$X_M:=\sum_e\chi_c(\Gr_e(M_0))\prod_i x_i^{(B_Te)_i-\inprod{S_i,M}+m_i}.$$

Define the map $\lambda_?$ : $\cc\rightarrow \Z^n$ to be given by $M\mapsto (\inprod{S_i,M}-m_i)$.

\begin{conjecture}
The cones generated by the images of \tilting objects under $\lambda_?$ provide a complete simplicial fan.

\end{conjecture}
Note that if we replace the map $\lambda_?$ by the dimension vector map,
we obtain a complete fan which is in general not simplicial.

\end{subsection}
\begin{subsection}{}

We  finish with a positivity conjecture which can be seen as an analogue of
Lusztig's positivity theorem, \cf \cite{lusztigbook}, for the dual canonical
basis.

\begin{conjecture}
For any object $M$ and any \exceptional object $K$ of $\cc$, the
integer $r(M,K)$ is non negative.
\end{conjecture}

In particular, the conjecture implies that the coefficients in Theorem \ref{maintheorem2}
are positive. Note that the rational numbers $r(M,N,\Gamma)$ defined in section \ref{filtration}
can be negative.

\end{subsection}

\section{Appendix on constructibility}{}\label{append}

We present a general proof for the constructibility of the sets
$\Ext^1(M,N)_Y$ in a triangulated category with finitely
many isoclasses of indecomposables.\smallskip

Let $k$ be a field and $\ct$ a $k$-linear triangulated category
with suspension functor $S$ such that
\begin{itemize}
\item[-] all $\Hom$-spaces in $\ct$ are finite-dimensional,
\item[-] each indecomposable of $\ct$ has endomorphism ring $k$,
\item[-] each object is a finite direct sum of indecomposables,
\item[-] $\ct$ has Serre duality, \ie there is an equivalence
$\nu: \ct\to\ct$ such that we have
\[
D\Hom(X,?) \iso \Hom(?,\nu X)
\]
for each $X\in\ct$, where $D$ denotes the functor $\Hom_k(?,k)$.
\end{itemize}
It is not hard to show that the last condition is a consequence
of the first three.
The conditions imply that $\ct$ has Auslander-Reiten
triangles and that the Auslander-Reiten translation $\tau$ is
given by $S^{-1} \nu$.

For objects $X,Y,Z$ of $\ct$, let $\Hom(X,Y)_Z$ be the set of
morphisms $f:X \to Y$ such that there is a triangle
\[
\xymatrix{ X \ar[r]^f & Y \ar[r] & Z \ar[r] & SX .}
\]

\begin{proposition} The set $\Hom(X,Y)_Z$ is constructible.
\end{proposition}

\begin{proof} Recall that a {\em split} triangle is a triangle which
is a direct sum of triangles one of whose three morphisms is an
isomorphism. Let us call  a triangle {\em minimal} if it does not
have a non zero split triangle as a direct factor. A triangle
\[
\xymatrix{ X \ar[r]^f & Y \ar[r] & Z \ar[r] & SX .}
\]
is minimal iff, in the category of morphisms, $f$ does
not admit non zero factors of the following forms
\[
\xymatrix{U \ar[r]^{\id_U} & U} \ko
\xymatrix{U \ar[r] & 0} \ko
\xymatrix{0 \ar[r] & U}.
\]
Let us call  such morphisms $f$ {\em minimal}. We now proceed by
induction on the sum $s$ of the numbers of indecomposable modules
occurring in the decompositions of $X$ and $Y$ into
indecomposables. Clearly the assertion holds if $s=0$ \ie $X=Y=0$.
Let us suppose $s>0$. Then $\Hom(X,Y)_Z$ is the disjoint union of
the set $M$ of morphisms $f$ such that the triangle
\[
\xymatrix{ X \ar[r]^f & Y \ar[r] & Z \ar[r] & SX .}
\]
is minimal and of the set $M'$ of morphisms
admitting a non zero direct factor of one of the above forms. The
set $M'$ is the union of orbits under $\Aut(X) \times \Aut(Y)$
of morphisms of the forms
\[
f' \oplus \id_{U}: X' \oplus U \to Y' \oplus U \ko
[f' , 0] : X' \oplus U \to Y' \ko
[f', 0]^t : X' \to Y' \oplus U \ko
\]
where $U$ is non zero and $f'$ runs through the sets $\Hom(X',
Y')_{Z'}$ for suitable $X'$, $Y'$, $Z'$, of which there are only a
finite number. It therefore follows from the induction hypothesis
that $M'$ is constructible. It remains to be shown that the set
$M$ is constructible. We work in the category $\mod \ct$ of
finitely presented functors on $\ct$ with values in the category
of $k$-vector spaces. It is an abelian category. Its projective
objects are the representable functors $\hat{U}=\Hom(?,U)$ and
these are also the injective objects. If $U$ is indecomposable in
$\ct$ and $S_U$ is the simple top of the indecomposable projective
$\hat{U}$, then $S_U$ admits the minimal projective presentation
\[
\xymatrix{ \hat{E_U} \ar[r]^{\hat{p_U}} & \hat{U} \ar[r] & S_U
\ar[r] & 0 \ko}
\]
where
\[
\xymatrix{ \tau U \ar[r] & E_U \ar[r]^{p_U} & U \ar[r] & S \tau U}
\]
is an Auslander-Reiten triangle. Moreover, $S_U$ is also the
simple socle of the indecomposable injective
\[
\hat{S\tau U} = \hat{\nu U}.
\]
If $f: X \to Y$ is a minimal morphism, then the morphisms
\[
\hat{X} \to \im(\hat{f}) \ko
\im(\hat{f}) \to \hat{Y}
\]
induced by $f$ are a projective cover and an injective hull, respectively.
Moreover, if $f$ is minimal and
\[
\xymatrix{X \ar[r]^f & Y \ar[r]^g & Z \ar[r] & SX}
\]
is a triangle, then $g$ is also minimal, so that $\hat{Z}$ is an
injective hull of $\im(\hat{g}) = \cok(\hat{f})$. Therefore, the
multiplicity $m_U$ of an indecomposable object $\nu U$ in the
decomposition of $Z$ into indecomposables equals the multiplicity
of the simple $S_U$ in the socle of $\cok(\hat{f})$. Since $U$ has
endomorphism algebra $k$, this also holds for $S_U$ and the
multiplicity $m_U$ equals
\[
\dim \Hom(S_U, \cok(\hat{f})).
\]
Now we have projective presentations
\[
\xymatrix{\hat{E_U} \ar[r]^{\hat{p_U}} &  \hat{U} \ar[r] &  S_U
\ar[r] & 0}
\]
and
\[
\xymatrix{\hat{X} \ar[r]^{\hat{f}} & \hat{Y} \ar[r] &
                                       \cok(\hat{f}) \ar[r] & 0.}
\]
Thus, we can compute the space $\Hom(S_U, \cok(\hat{f}))$ as the
quotient of the space of morphisms
\[
\xymatrix{E_U \ar[r]^{p_U} \ar[d]_a & U  \ar[d]^b  \\
X \ar[r]^f & Y}
\]
modulo the subspace formed by the morphisms of the
form $(a,b) = (cp_U, fc)$ for some morphism $c: U \to X$. The condition
\[
\dim \Hom(S_U, \cok(\hat{f})) = m_U
\]
then translates into conditions on the ranks of the linear maps
\[
(a,b) \mapsto bp_U - fa \mbox{ and } c \mapsto (cp_U, fc).
\]
Clearly, the $f\in \Hom(X,Y)$ satisfying these rank conditions
for all indecomposables $U$ form a constructible subset.
The intersection of this subset with the complement of $M'$
is still constructible (since $M'$ is) and clearly equals $M$.
So $M$ is constructible.
\end{proof}

\end{subsection}

\end{document}